\documentclass[11pt,amssymb]{amsart}
\usepackage{amsmath,amssymb,url,hyperref,amsthm}
\usepackage[mathscr]{eucal}
\usepackage[all,cmtip]{xy}
\usepackage{tikz-cd}
\usepackage{enumitem}

\newcommand{\im}{{\rm im}\:}

\newcommand{\Z}{{\mathbb Z}}
\newcommand{\cJ}{\mathcal{J}}

\newcommand{\Q}{{\mathbb Q}}

\newcommand{\Nrd}{\mathrm{Nrd}}

\newcommand{\R}{{\mathbb R}}

\newcommand{\Br}{\mathrm{Br}}

\newcommand{\cO}{\mathscr{O}}

\newcommand{\Ga}{\mathrm{Gal}}
\newcommand{\GL}{\operatorname{GL}}
\newtheorem{thm}{Theorem}[section]
\newtheorem{lemma}[thm]{Lemma}
\newtheorem{prop}[thm]{Proposition}
\newtheorem{cor}[thm]{Corollary}
\theoremstyle{definition}
\newtheorem{definition}[thm]{Definition}
\newtheorem{remark}[thm]{Remark}
\newtheorem*{remark*}{Remark}

\newcommand{\gen}{\mathbf{gen}}
\newcommand{\Hom}{\mathrm{Hom}}

\newcommand{\Pic}{\mathrm{Pic}}
\newcommand{\Spec}{\operatorname{Spec}}
\newcommand{\cQ}{\mathcal{Q}}
\newcommand{\tv}{\tilde v}
\newcommand{\he}{H_{\mathrm{\acute{e}t}}}
\newcommand{\uG}{\underline{G}}

.240pk scaled 1200 .240pk
\DeclareFontFamily{U}{wncy}{}
    \DeclareFontShape{U}{wncy}{m}{n}{<->wncyr10}{}
    \DeclareSymbolFont{mcy}{U}{wncy}{m}{n}
    \DeclareMathSymbol{\Sha}{\mathord}{mcy}{"58}
\begin{document}

\title[Groups with good reduction, genus problem, and unramified cohomology]{Simple algebraic groups with good reduction, the genus problem, and unramified cohomology}

\author[I.A.~Rapinchuk]{Igor A. Rapinchuk}

\begin{abstract}
We obtain a number of new results on the arithmetic of linear algebraic groups over higher-dimensional fields. In particular, we resolve some key open questions over 2-dimensional global fields, such as the finiteness of the genus of spinor groups of quadratic forms of dimension $\geq 10$. We also extend various finiteness results to new classes of fields, including function fields of $p$-adic surfaces. Our arguments rely in a crucial way on new finiteness results for unramified cohomology.
\end{abstract}

\address{Department of Mathematics, Michigan State University, East Lansing, MI
48824, USA}

\email{rapinchu@msu.edu}

\maketitle

\tableofcontents

\section{Introduction}\label{S:Intro}

This paper is a new installment in the study of the arithmetic of linear algebraic groups over fields more general than local and global fields, with a special focus on finiteness properties. Our principal contributions include the following. On the one hand, we resolve some key issues over 2-dimensional global fields that were left open in \cite{Spinor}, such as the finiteness of the genus of spinor groups of quadratic forms of dimension $\geq 10$ --- see Theorem \ref{T-GenusSpinor}. On the other hand, we extend a number of finiteness results to new classes of fields, including function fields of $p$-adic surfaces (cf. Theorem \ref{T-Surface}). Both of these objectives are achieved mainly via new finiteness results for unramified cohomology (see Theorem \ref{T-SBUnram} and Proposition \ref{P-FinUnramPadic}).


One of the central elements in the current work on the arithmetic of linear algebraic groups is the analysis of {\it groups with good reduction} and associated Shafarevich-type finiteness conjectures (see \cite{RR-survey} for an overview). Thus, broadly speaking, one would like to develop a general understanding of the (connected, reductive) linear algebraic groups over a field $K$ that have good reduction at a suitable set $V$ of discrete valuations of $K$, particularly with the goal of addressing finiteness problems for the isomorphism classes of $K$-forms with good reduction of a given reductive group. More specifically, it is expected that the number of such forms is finite in the following two situations:

\vskip3mm

\noindent (AFC) \parbox[t]{15cm}{$K$ is a finitely generated field (of ``good" characteristic\footnotemark) and $V$ is a divisorial set of places of $K$ (Arithmetic Finiteness Conjecture, cf. \cite[Conjecture 5.7]{RR-survey})}

\footnotetext{In this context, we say that $p = \mathrm{char}\: K > 0$ is good for a semisimple $K$-group $G$ if $p$ does not divide the order of the Weyl group of $G$. For reductive non-semisimple groups, only characteristic zero is good.}

\vskip3mm

\noindent (GFC) \parbox[t]{15cm}{$K = k(X)$ is the function field of a geometrically integral normal variety $X$ defined over a field $k$ that has characteristic 0 and is of type (F) as defined by Serre (cf. \cite[Ch. III, \S 4]{Serre-GC}) and $V$ is the set of discrete valuations of $K$ associated with the prime divisors on $X$ (Geometric Finiteness Conjecture).}

\vskip3mm

Numerous results on these conjectures have been obtained in \cite{Spinor}, \cite{CRR-MaxTori}, \cite{RR1}, \cite{RR-Tori}, and \cite{RR3}. Here we only mention that both conjectures have been settled for inner forms of type $\mathsf{A}_{n-1}$ provided that $n$ is prime to the characteristic. Furthermore, in \cite{Spinor}, we proved (AFC) for spinor groups of quadratic forms over {\it 2-dimensional global fields}\footnote{Following Kato \cite{Kato}, by a 2-dimensional global field, we mean either the function field $k(C)$ of a smooth geometrically integral curve $C$ over a number field $k$, or the function field $k(S)$ of a smooth geometrically integral surface $S$ over a finite field $k$.} of characteristic $\neq 2.$ While this result can be extended to the special unitary groups of {\it hermitian} forms associated with quadratic extensions and quaternion algebras (equipped with the canonical involution), the case of {\it skew-hermitian} forms over quaternions remained open. Our first goal in the present paper is to resolve this case.



\begin{thm}\label{T-SB}
Let $K = k(C)$ be a 2-dimensional global field of characteristic 0 and suppose $V$ is a divisorial set of places of $K$. Then the set of $K$-isomorphism classes of the universal covers $\tilde{G}$ of the special unitary groups $G = \mathrm{SU}_n(h,D)$, where $h$ is a nondegenerate $n$-dimensional skew-hermitian form over some quaternion division $K$-algebra $D$ with canonical involution,
having good reduction at all $v \in V$ is finite.
\end{thm}

(We should point out that (GFC) for such groups defined over function fields of curves over fields of type (F) was proved in \cite{Srimathy}.)

\vskip2mm

Next, let us describe two important applications of this result. In \cite{Weis}, Weisfeiler initiated the study of absolutely almost simple algebraic groups that are defined over a field $K$ and become split over a (separable) quadratic extension $L/K$. This class of groups, on the one hand, generalizes the classical case  $K = \R$ that goes back to Cartan, and on the other hand, is in some sense ``closest'' to split groups over general fields. For these reasons, it provides a natural ``testing ground" for various conjectures. While some structural results for such groups have been obtained (cf. \cite{Chernousov}, \cite[\S9.4]{PR}, \cite{RapCSP}), their complete classification (for example, in terms of cohomological or Witt ring invariants) is known only over the real numbers (see, e.g., \cite{Borovoi} for a modern treatment). Developing such a classification would be a task of considerable independent interest, which would also be helpful in the analysis of the above finiteness conjectures in the general case. However, here is what has been accomplished so far by {\it ad hoc} considerations. The results of \cite{Spinor} imply (AFC) over 2-dimensional global fields for absolutely almost simple algebraic groups splitting over a quadratic extension of the base field of types $\textsf{A}_{\ell}$, $\textsf{B}_{\ell}$, $\textsf{C}_{\ell}$, and $\textsf{G}_2$; recently, in \cite{CRR-MaxTori}, we considered type $\textsf{F}_4$. Theorem \ref{T-SB}
yields an extension of this result to groups of type $\textsf{D}_{\ell}$ splitting over a quadratic extension (we recall that a simply connected group of type $\textsf{D}_{\ell}$ that splits over a quadratic extension is either the spinor group of a quadratic form or the universal cover of a special unitary group as in Theorem \ref{T-SB} --- see \cite[Ch.~2]{PRR}),
leaving open only types $\textsf{E}_{\ell}$, for $\ell = 6, 7$, and $8$.

Another application of Theorem \ref{T-SB} is to the genus problem. We recall that
the genus $\gen_K(G)$ of an absolutely almost simple simply connected algebraic $K$-group $G$ is defined as the set of $K$-isomorphism classes of {inner} $K$-forms of $G$ that have the same isomorphism classes of maximal $K$-tori as $G$ (see \cite[\S1]{CRR-MaxTori} and \cite[\S8]{RR-survey} for further details). It has been conjectured
(see \cite[Conjecture 8.10]{RR-survey}) that $\gen_K(G)$ should be finite whenever the base field $K$ is finitely generated (of ``good'' characteristic). In \cite[Theorem 1.1]{CRR-MaxTori}, we linked this conjecture with (AFC); however, the proof of the latter for spinor groups of quadratic forms over 2-dimensional global fields given in \cite{Spinor}
was insufficient to establish the finiteness of the genus of such groups. Combining the results of \cite{Spinor} and \cite{CRR-MaxTori} with Theorem \ref{T-SB}, we are now able to provide a complete finiteness statement for spinor groups of quadratic forms of sufficiently large dimension.

\begin{thm}\label{T-GenusSpinor}
Suppose $K = k(C)$ is a 2-dimensional global field of characteristic 0, and let $q$ be a nondegenerate quadratic form over $K$ in $n \geq 10$ variables. Then the genus $\gen_K(G)$ of the spinor group $G = {\rm Spin}_n(q)$ is finite.
\end{thm}

We will now give some indications of the proof of Theorem \ref{T-SB}, highlighting the connection with unramified cohomology. Suppose $D$ is a quaternion division algebra over $K = k(C)$ and let $Q$ be the corresponding Severi-Brauer conic. Since the function field $L = K(Q)$ is a splitting field of $D$, by Morita equivalence, to any skew-hermitian form $h$ over $D$, one can associate a quadratic form $q_{h_L}$ over $L$ (see \S\ref{S:Morita}). Then one relates the study of the universal covers $\tilde{G} = \widetilde{{\rm SU}_n(h,D)}$ over $K$ having good reduction to the analysis of quadratic forms over $L$ with good reduction. Using the Milnor Conjecture (proved by Voevodsky) and some idelic considerations, one shows that the finiteness properties of the latter are intimately connected with the finiteness properties of the unramified cohomology of $L$ with $\mu_2$-coefficients (see Proposition \ref{P:FormsGoodReduction}). So, to conclude the argument, the crucial ingredient is the following result.

\begin{thm}\label{T-SBUnram}
Let $K = k(C)$ be a 2-dimensional global field of characteristic 0. Suppose $$Q: ax^2 + by^2 = z^2$$ (with $a, b \in K^{\times}$) is an anisotropic conic in $\mathbb{P}^2_K$,
and
consider the function field $L = K(Q).$ Then for any divisorial set of places $V$ of $L$, the
unramified cohomology groups $H^i(L, \mu_2)_{V}$ are finite for all $i \geq 1.$
\end{thm}

This theorem can be viewed as a new contribution to the notoriously difficult and largely open problem of studying the finiteness properties of the degree 3 unramified cohomology of surfaces over number fields. As we will see in \S\ref{S:FinUnram}, the
most involved part of the proof is showing the finiteness of $H^3(L, \mu_2)_{{V}}$, for which we use, among other considerations, the Bloch-Kato Conjecture and the calculation of the unramified cohomology of quadrics by Kahn, Rost, and Sujatha \cite{KRS}.


Turning our attention now to the geometric case, we obtain the following result concerning (GFC) (by a $p$-adic field, we mean a finite extension of $\Q_p$).

\begin{thm}\label{T-Surface}
Suppose that
$p$ is an odd prime and
$S$ is a smooth geometrically integral surface over a $p$-adic field $k$.
Let $K = k(S)$ be the function field of $S$ and $V$ be the set of geometric places of $K$ corresponding to all of the prime divisors of $S$. Then we have the following.

\vskip1mm

\begin{enumerate}[label=\textup{(\roman*)}, leftmargin=*]
  \item Let $n \geq 5$ be an integer. Then the set of $K$-isomorphism classes of spinor groups $G = \mathrm{Spin}_n(q)$ of nondegenerate quadratic forms in $n$ variables over $K$ that have good reduction at all $v \in V$ is finite.
  \item Let $L/K$ be a quadratic extension and $n \geq 2$ be an integer. Then the set of $K$-isomorphism classes of special unitary groups $G = {\rm SU}_n(L/K, h)$ of nondegenerate hermitian $L/K$-forms $h$ in $n$ variables that have good reduction at all $v \in V$ is finite.
  \item Let $G$ be a simple algebraic $K$-group of type $\mathsf{G}_2.$ Then the set of $K$-isomorphism classes of $K$-forms $G'$ of $G$ that have good reduction at all $v \in V$ is finite.
  \item The set of $K$-isomorphism classes of $K$-forms of type $\mathsf{F}_4$ that split over a quadratic extension of $K$ and have good reduction at all $v \in V$ is finite.
\end{enumerate}

\end{thm}

Whereas these statements were previously known when $K = k(C)$ is the function field of a {\it curve} (cf. \cite{Spinor}, \cite{CRR-MaxTori}), the significance of Theorem \ref{T-Surface} is that it takes one of the first steps in investigating these kinds of finiteness problems over the function fields of {\it higher-dimensional} varieties. The proof of the theorem relies on the finiteness of the unramified cohomology groups $H^i(k(S), \mu_2)_V$ $(i \geq 1$),
for which, at present, we need to restrict ourselves to the case where $k$ is a $p$-adic field (the assumption that $p$ is odd is in fact needed only for the finiteness of $H^3(k(S), \mu_2)_V$). In future work, we plan to address the case where $k$ is any field that has characteristic 0 and is of type (F).

As we already remarked above, our results in \cite{CRR-MaxTori} establish a link between (AFC) and the genus problem. In order to apply a similar strategy to the study of the genus problem in the geometric case, an important step would be to extend \cite[Theorem 1.1]{CRR-MaxTori} to the situation where the relevant residue field is not necessarily finitely generated.
In this direction, we prove the following.


\begin{thm}\label{T-GenusHilbert}
Suppose $X$ is a smooth geometrically integral variety of dimension $> 1$ over a field $k$ of characteristic 0. Let $K = k(X)$ be the function field of $X$ and $v$ be a geometric place of $K$. If $G$ is an absolutely almost simple algebraic group over $K$ that has good reduction at $v$, then every $G' \in \gen_K(G)$ also has good reduction at $v$.

\end{thm}

We remark that, although the proof of this statement is inspired by the techniques developed in \cite{CRR-MaxTori}, the arguments require new considerations involving generic tori over Hilbertian fields, which are of independent interest. 


With Theorem \ref{T-GenusHilbert} at our disposal, we are then able to apply Theorem \ref{T-Surface}, in conjunction with some considerations from \cite{Spinor}, to obtain several new statements concerning the genus problem for groups over the function fields of $p$-adic surfaces.

\begin{thm}\label{T-GenusSurface}
Let $p$ be an odd prime and $K = k(S)$ be the function field of a smooth geometrically integral surface $S$ over a $p$-adic field $k$. Then $\gen_K(G)$ is finite in the following cases:

\vskip1mm

\begin{enumerate}[label=\textup{(\roman*)}, leftmargin=*]
  \item $G = {\rm Spin}_n(q)$, where $q$ is an $n$-dimensional nondegenerate quadratic form over $K$ and either $n \geq 5$ is odd or $n \geq 10$ is even and $q$ is $K$-isotropic;
  \item $G$ is a simple algebraic $K$-group of type $\mathsf{G}_2$;
  \item $G$ is a simple algebraic $K$-group of type $\mathsf{F}_4$ that splits over a quadratic extension of $K$.
\end{enumerate}

\end{thm}

\vskip2mm

The final topic that we consider in this paper is (higher) cohomological
local-global principles for algebraic tori. We recall that if $K$ is a field equipped with a set $V$ of discrete valuations and $G$ is a linear algebraic group defined over $K$, then one considers the global-to-local map in Galois cohomology
$$
\lambda_{G,V} \colon H^1(K, G) \to \prod_{v \in V} H^1(K_v, G)
$$
and defines the corresponding Tate-Shafarevich set as
$$
\Sha^1(K, G, V) = \ker \lambda_{G,V}.
$$
One of the important problems in the arithmetic theory of linear algebraic groups is to determine when $\lambda_{G,V}$ is injective (and hence $\Sha^1(K,G,V)$ is trivial), or at least when it is proper (which would imply that $\Sha^1(K,G,V)$ is finite). In the classical setting of number fields (with $V$ being the set of all places, including the archimedean ones), perhaps the most significant results are that $\Sha^1(K,G,V)$ is finite when $G$ is a torus and it is trivial whenever $G$ is a semisimple simply connected or adjoint group. More generally, the expectation is that $\lambda_{G,V}$ is proper for {\it any reductive group} in both the arithmetic case of finitely generated fields and the geometric case of function fields of varieties over fields of type (F) --- we will call this the Properness Conjecture.
We refer the reader to \cite[\S6]{RR-survey} for an overview of cohomological local-global principles as well as a discussion of the connections between the Properness Conjecture and groups with good reduction, and to \cite{Spinor}, \cite{RR-Tori}, and \cite{RR3} for recent results in this direction. Here, we note that $\lambda_{G,V}$ is indeed proper when $G$ is a torus in both the arithmetic and geometric cases. Our purpose in the present work is to establish the finiteness of {\it higher} Tate-Shafarevich groups of algebraic tori in the geometric setting.

Before formulating our result, let us point out that higher Tate-Shafarevich groups arise naturally in arithmetic contexts. Indeed, if $L/K$ is a finite Galois extension of number fields and $V^K$ is the set of all places of $K$, then the injectivity of the global-to-local map
$$
\beta_{L/K, V^K} \colon H^2(L/K, L^{\times}) \longrightarrow \prod_{v \in V^K} H^2 (L_{\bar{v}}/K_v, L_{\bar{v}}^{\times})
$$
(where for each $v \in V^K$, we fix a \emph{single} extension $\bar{v}$ of $v$ to $L$) is the famous theorem of Albert-Brauer-Hasse-Noether, which is essentially equivalent to the Reciprocity Law in Global Class Field Theory (cf. \cite[Ch. VII]{ANT}). From this result, one then obtains the finiteness of $\ker \beta_{L/K, V}$ for any set $V$ that contains almost all places of $K$.

Our work shows that the higher Tate-Shafarevich groups of tori are finite in much greater generality. To fix notations, suppose $T$ is a torus over a field $K$ that is equipped with a set $V$ of discrete valuations. Let $L/K$ be a finite Galois extension and denote by $\bar{V}$ the set of all extensions of valuations from $V$ to $L.$ Then for all $i \geq 1$, we define
$$
\Sha^i(L/K, T, V) = \ker \left(H^i(L/K, T) \to \prod_{v \in V} H^i(L_{\bar{v}}/K_v, T) \right),
$$
where, as above, for each $v \in V$, we choose a single $\bar{v} | v$ in $\bar{V}$ (we note that when $i =1$ and $L/K$ is a finite Galois splitting field of $T$, then $\Sha^1(L/K, T,V) = \Sha^1(K,T,V)$ by Shapiro's Lemma and Hilbert's Theorem 90).

In \cite[Proposition 4.2]{RR-Tori}, we considered the arithmetic situation and showed that when $V$ is a divisorial set of places of a finitely generated field $K$, then $\Sha^i(L/K, T, V)$ is finite for any finite Galois extension $L/K$, any $K$-torus $T$, and all $i \geq 1$. Then, in \cite[Theorem 1.5]{RR3}, we analyzed the geometric case where $K = k(X)$ is the function field of a normal irreducible
variety $X$ over a field $k$ of characteristic zero and of type (F) and $V$ is a set of geometric places of $K$, and established the finiteness of $\Sha^1(K,T,V)$. We now complete this line of work and establish the finiteness of $\Sha^i$ for all $i \geq 2$ in the geometric situation.

\begin{thm}\label{T-HigherSha}
Suppose $k$ is a field of characteristic 0 and of type $\mathrm{(F)}$, let $K = k(X)$ be the function field of a smooth geometrically integral variety $X$ over $k$, and denote by $V$ the set of geometric places of $K$ associated with the prime divisors of $X$. Let $T$ be a torus over $K$ and $L/K$ be a finite Galois extension that splits $T$. Then the higher Tate-Shafarevich groups $\Sha^i (L/K, T, V)$ are finite for all $i \geq 2.$

\end{thm}

\vskip3mm

\noindent {\bf Notations and conventions.} For a discrete valuation $v$ of a field $F$, we will denote by
$F_v$ and $F^{(v)}$ the associated completion and residue field, respectively. Furthermore, we let $\mathcal{O}_{F,v} \subset F$ and
$\mathcal{O}_{F_v} \subset F_v$ denote the corresponding valuation rings.

An affine reductive algebraic group $G$ over $F$ is said to have \emph{good reduction} at $v$ if there exists a reductive group scheme $\mathscr{G}$ over $\mathcal{O}_{F_v}$ with generic fiber $\mathscr{G}  \times_{\mathcal{O}_{F_v}} F_v$ isomorphic to $G \times_F F_v$. The closed fiber $\mathscr{G} \times_{\mathcal{O}_{F_v}} F^{(v)}$ is called the {\it reduction} of $G$ at $v$ and will be denoted $\uG^{(v)}$ (we note that by the results of \cite{Guo} and \cite{Nisn2}, the reductive model $\mathscr{G}$ is unique, and hence the reduction is well-defined). It is known that in the case where $G$ is an absolutely almost simple $F$-group, the reduction $\uG^{(v)}$ has the same Killing-Cartan type as $G$ (see \cite[Exp. XXII, Proposition 2.8]{SGA3}).

We will follow the usual conventions concerning Galois cohomology (cf. \cite{Serre-GC}). In particular, given an algebraic group $G$ defined over a field $F$ and a Galois extension $L/F$, we denote by $H^1(L/F, G)$ the pointed set $H^1(\Ga(L/F), G(L))$ of (continuous) noncommutative Galois cohomology, and we write $H^1(F, G)$ for $H^1(F^{\rm sep}/F, G(F^{\rm sep}))$, where $F^{\rm sep}$ is a separable closure of $F$. In the case of commutative coefficients, given a $\Ga(L/F)$-module $M$ (assumed to be discrete if $L/F$ is infinite), we will write $H^i(L/F, M)$ for the Galois cohomology group $H^i(\Ga(L/F), M)$;
when $L = F^{\rm sep}$, we simply write $H^i(F, M).$ Moreover, if $T$ is an $F$-torus, we  denote by $H^i(L/F, T)$ the Galois cohomology group $H^i(\Ga(L/F), T(L))$, which we again abbreviate to $H^i(F, T)$ when $L = F^{\rm sep}.$
For extra clarity, we will occasionally write $H^i(L/F, T(L)).$

Finally, throughout the paper, we will use $F$ to denote a general field (often equipped with a discrete valuation or a set of discrete valuations), whereas $K$ will be reserved for the special classes of fields under consideration, namely finitely generated fields (including 2-dimensional global fields) and function fields of varieties over fields of type (F) (including function fields of surfaces over $p$-adic fields).

\section{Some basic facts about divisorial sets, adeles, and unramified cohomology}\label{S:Recollections}


The purpose of this section is to give an overview of some background material that will be used throughout the rest of the paper.

\vskip2mm

\noindent {\bf 2.1. Divisorial sets of places.} Let $K$ be a field that is finitely generated over its prime subfield (in the sequel, we refer to such fields simply as {\it finitely generated}). It is well known that $K$ admits a model $\mathcal{X}$, i.e. an irreducible separated normal scheme over $\Z$ with function field $K$.
One then considers the set $V = V_{\mathcal{X}}$ of discrete valuations of $K$ corresponding to all the prime divisors on $\mathcal{X}$. Any set of discrete valuations of $K$ having this form will be called {\it divisorial}.
In more algebraic terms, we observe that finding a model of $K$ essentially amounts to specifying an integrally closed subring $A \subset K$ that is a finitely generated $\Z$-algebra; we then have the set of discrete valuations of $K$ associated with the height one prime ideals of $A$. We note that given {\it any model} $\mathcal{X}$ of $K$ with corresponding divisorial set $V_{\mathcal{X}}$, one can consider an affine open subscheme $\mathcal{U} = \mathrm{Spec}\: A$ of $\mathcal{X}$. Then for the divisorial set $V' = V_{\mathcal{U}}$ associated with $\mathcal{U}$, we have $V' \subset V$ and $V \setminus V'$ is finite. So, for our purposes, we can generally restrict our attention to divisorial sets associated with {\it affine} models. Two important properties of divisorial sets are the following (cf. \cite[\S5.3]{RR-survey}):

\vskip1mm

\noindent $\bullet$ \parbox[t]{16.5cm}{Any two divisorial sets $V_1$ and $V_2$ associated with two different models of $K$ are commensurable, i.e. $V_i \setminus (V_1 \cap V_2)$ is finite for $i = 1,2$.}

\vskip1mm

\noindent $\bullet$ \parbox[t]{16.5cm}{For any finite subset $S$ of a divisorial set $V$, the set $V \setminus S$ contains a divisorial set. More precisely, if $V$ is associated with a model $\mathcal{X}$, then we can find an open subscheme $\mathcal{U} \subset \mathcal{X}$ (which is simply a suitable principal open set $\mathcal{U} = {\rm Spec}~A_f$ in the case where $\mathcal{X} = {\rm Spec}~A$ is affine) such that
for the corresponding divisorial set $V_{\mathcal{U}}$, we have $V_{\mathcal{U}} \subset V \setminus S$.}

\vskip3mm

\noindent {\bf 2.2. Adeles.} Suppose $F$ is a field equipped with a set $V$ of discrete valuations satisfying
the following condition:
\vskip2mm

\noindent (*) for any $a \in F^{\times}$, the set $V(a) := \{ v \in V \mid v(a) \neq 0 \}$ is finite.

\vskip2mm

\noindent (For example, we note that (*) automatically holds for any divisorial set $V$ of a finitely generated field $K$.)

\vskip1mm
We denote by $\mathbb{A}(F , V)$ the corresponding ring of adeles, defined as the restricted product of the completions $F_v$ for $v \in V$ with respect to the valuation rings $\mathcal{O}_{F_v} \subset F_v$:
$$
\mathbb{A}(F,V) = \left\{ (x_v) \in \prod_{v \in V} F_v \mid x_v \in \mathcal{O}_{F_v} \ \text{for almost all} \ v \in V \right\}.
$$
The group of ideles $\mathbb{I}(F,V)$ is by definition the group of invertible elements of $\mathbb{A}(F,V)$.
In view of (*), we have diagonal embeddings $F \hookrightarrow \mathbb{A}(F,V)$ and $F^{\times} \hookrightarrow \mathbb{I}(F,V).$ In particular, $\mathbb{A}(F,V)$ has a natural structure of an $F$-algebra.
Furthermore, we let
$$
\mathbb{A}^{\infty}(F,V) = \prod_{v \in V} \mathcal{O}_{F_v}  \ \ \ \text{and} \ \ \ \mathbb{I}^{\infty}(F,V) = \prod_{v \in V} \mathcal{O}_{F_v}^{\times}
$$
denote the subring of integral adeles and subgroup of integral ideles, respectively. We note that these constructions behave in the expected way
with respect to finite Galois extensions. Namely, suppose $L/F$ is a finite Galois extension with Galois group $\mathcal{G} = \Ga(L/F)$, and let $\bar{V}$
denote the set of all extensions to $L$ of the discrete valuations in $V$. 
It is well-known that we have a natural isomorphism
$$
\mathbb{A}(K,V) \otimes_K L \simeq \mathbb{A}(L, \bar{V})
$$
of topological rings (see, e.g., \cite[Ch II, \S 14]{ANT} for the number field case), which, in particular, enables us to define an action of $\mathcal{G}$ on $\mathbb{A}(L, \bar{V})$ with the property that
$$
\mathbb{A}(L, \bar{V})^{\mathcal{G}} = \mathbb{A}(K,V).
$$

Now, given an $F$-torus $T$, since $\mathbb{A}(F,V)$ is an $F$-algebra, we can consider the group of adelic points $T(\mathbb{A}(F,V)).$ It turns out that this group may be described more concretely as follows. For each $v \in V$, let $T(\mathcal{O}_{F_v})$ be the unique maximal bounded subgroup of $T(F_v)$. Then one can show that $T(\mathbb{A}(F,V))$ is the restricted product of the groups $T(F_v)$ for $v \in V$ with respect to the subgroups $T(\mathcal{O}_{F_v})$ (see \cite[\S11.1]{Voskr} for the argument in the number field case). We then define the subgroup of integral adeles to be
$$
T(\mathbb{A}^{\infty}(F,V)) = \prod_{v \in V} T(\mathcal{O}_{F_v}).
$$
We note that the diagonal embedding $F \hookrightarrow \mathbb{A}(F,V)$ yields an embedding $T(F) \hookrightarrow T(\mathbb{A}(F,V))$, which one usually also calls the diagonal embedding. For a finite Galois extension $L/F$, the corresponding groups $T(\mathbb{A}(L, \bar{V}))$ and $T(\mathbb{A}^{\infty}(L, \bar{V}))$ are defined analogously, and, moreover, the action of $\mathcal{G} = \Ga(L/F)$ on $\mathbb{A}(L, \bar{V})$ induces a $\mathcal{G}$-action on $T(\mathbb{A}(L, \bar{V})).$ In particular, the diagonal embedding gives rise to group homomorphisms
\begin{equation}\label{E:ShaAdeleHom}
\nu^i_{L/F} \colon H^i(L/F, T) \to H^i(L/F, T(\mathbb{A}(L, \bar{V})))
\end{equation}
for all $i \geq 1.$

\vskip3mm

\noindent {\bf 2.3. Picard groups.} As above, suppose $F$ is a field equipped with a set $V$ of discrete valuations satisfying condition (*).
Set ${\rm Div}(F,V)$ to be the free abelian group on $V$ (``group of divisors"), and let $\mathrm{Pic}(F , V)$ be the {\it Picard group} of $V$, defined as the quotient of ${\rm Div}(F,V)$ by the subgroup of principal divisors $\mathrm{P}(F , V)$, where the principal divisor corresponding to $a \in F^{\times}$ is $\sum_{v \in V} v(a) \cdot v$ (observe that this definition makes sense due to (*)). It is straightforward to see that the map
$$
\phi \colon \mathbb{I}(F,V) \to {\rm Div}(F,V), \ \ \ (x_v) \mapsto \sum_{v \in V} v(x_v) \cdot v
$$
induces an identification
\begin{equation}\label{E:PicIdeles}
\mathbb{I}(F,V)/ (\mathbb{I}^{\infty}(F,V) \cdot F^{\times}) \stackrel{\sim}{\longrightarrow} \Pic(F,V)
\end{equation}
(cf. \cite[Lemma 2.2]{Spinor}).


Suppose now that $K$ is a finitely generated field equipped with a divisorial set $V$ of discrete valuations associated with a model $\mathcal{X}$. As remarked above, after possibly deleting finitely many places from $V$, we may assume that $\mathcal{X} = {\rm Spec}~A_0$ is affine and the elements of $V$ correspond to the height one prime ideals of $A_0$.
It follows from \cite[Corollaire 1]{Kahn} that $\mathrm{Pic}(K,V)$ is finitely generated. Consequently, there exists a subset $V' \subset V$ with finite complement such that $\mathrm{Pic}(K , V')$ is trivial (see, for example, \cite[Proposition 3.1]{RR-Tori}). On the other hand,
as noted above, the cofinite set $V'$ contains a divisorial set $V_U$ corresponding to a principal open subset $U = {\rm Spec}~(A_0)_f$ (with $V' \setminus V_U$ a finite set). Then we also have $\Pic(K, V_U) = \{0 \}.$ From this, we see that since $(A_0)_f$ is noetherian and normal and ${\rm Cl}~(A_0)_f = \Pic(K, V_U) = \{0 \}$, the ring $(A_0)_f$ is in fact a unique factorization domain by \cite[Proposition 6.2]{Hartshorne}. Moreover,
in view of (\ref{E:PicIdeles}), the triviality of $\Pic(K, V_U)$ yields the equality
$$
\mathbb{I}(K,V_U) = \mathbb{I}^{\infty}(K,V_U) \cdot K^{\times}.
$$
For ease of reference, we record this discussion in the following statement.
\begin{lemma}\label{L:AdeleRecall}
Suppose $K$ is a finitely generated field equipped with a divisorial set $V$ of discrete valuations associated with a model $\mathcal{X}$. Then there exists a subset $V' \subset V$ with finite complement such that 
$$
\mathbb{I}(K,V') = \mathbb{I}^{\infty}(K,V') \cdot K^{\times}.
$$
In fact, for such $V'$, one can take the set of discrete valuations corresponding to the prime divisors  
of an affine open subset ${\rm Spec}~A \subset \mathcal{X}$, where $A$ is a noetherian unique factorization domain. 
\end{lemma}

\vskip3mm

\noindent {\bf 2.4. Vertical valuations and their properties.}
Let $F$ be a field equipped with a discrete valuation $v$, denote by $\mathcal{O}_{F,v} \subset F$ the corresponding valuation ring with maximal ideal $\mathfrak{m}$, fix a uniformizer $\pi \in \mathfrak{m}$, and set $F^{(v)} = \mathcal{O}_{F,v}/(\pi)$ to be the residue field. Let $\mathcal{X} \stackrel{\varphi}{\longrightarrow} {\rm Spec}~\mathcal{O}_{F,v}$ be a scheme that is of finite type and flat over $\mathcal{O}_{F,v}$, and suppose that $\mathcal{X}$ is normal and integral. We will assume that the generic fiber $\mathcal{X}_F$ is a smooth geometrically integral $F$-variety of dimension $n \geq 1$ and also that the special fiber $\mathcal{X}_{F^{(v)}} = \mathcal{X} \times_{\mathcal{O}_{F,v}} F^{(v)}$ is a smooth geometrically integral $F^{(v)}$-variety of dimension $n$. Finally, we denote by
$$
L = F(\mathcal{X}) = F(\mathcal{X}_F)
$$
the function field of $\mathcal{X}$.

Now, since $\mathcal{X}$ is flat over $\mathcal{O}_{F,v}$, a well-known consequence of Krull's principal ideal theorem is that the special fiber $\mathcal{X}_{F^{(v)}}$ has codimension 1 in $\mathcal{X}$ (see, e.g., \cite[Chapter 4, Theorem 3.12]{Liu}). Hence, if $\eta$ is the (unique) generic point of $\mathcal{X}_{F^{(v)}}$, then the local ring $\mathcal{O}_{\mathcal{X}, \eta}$ is a noetherian normal integral domain of dimension 1 and therefore is a discrete valuation ring with fraction field $L$ (cf. \cite[Proposition 9.2]{At}). We thus obtain a discrete valuation $\tilde{v}$ on $L$ that extends $v$, which is sometimes referred to as the {\it vertical valuation} of $\mathcal{X}$ over $v$. As usual, we define the ramification index $e(\tilde{v}|v)$ to be the index $[\tilde{v}(L^{\times}) : v(F^{\times})]$; equivalently, if we normalize $\tilde{v}$ to have value group $\Z$ (with $\tilde{v}(t) = 1$ for any uniformizer $t \in \mathcal{O}_{\mathcal{X}, \eta}$), then $e(\tilde{v}|v) = \tilde{v}(\pi).$

\begin{lemma}[{cf.\ \cite[Proposition 3.1(b)]{Harb-Stev}}]\label{lem22}
With notations as above, we have the following:
\begin{enumerate}[label=\textup{(\roman*)}, leftmargin=*]
  \item the ramification index is $e(\tv|v) = 1$;
  \item the residue field $L^{(\tv)}$ is the function field $F^{(v)}(\mathcal{X}_{F^{(v)}})$ of the special fiber;
  \item the extension $L^{(\tv)}/F^{(v)}$ is separably generated, i.e., it admits a transcendence basis $t_1, \dots, t_n$ such that $L^{(\tv)}$ is a finite separable extension of $F^{(v)}(t_1, \dots, t_n)$.
\end{enumerate}
\end{lemma}

\begin{proof} \ (i), (ii) Let $\Spec R \subset \mathcal{X}$ be an affine open subset that meets the special fiber $\mathcal{X}_{F^{(v)}}$. The ring $R/\pi R$ is the coordinate ring of $\Spec R \cap \mathcal{X}_{F^{(v)}}$, which is a nonempty open subset of the integral scheme $\mathcal{X}_{F^{(v)}}$. Hence $R/\pi R$ is an integral domain, so $\pi R$ is a prime ideal, and the corresponding point of $\Spec R$ is the generic point $\eta$ of $\mathcal{X}_{F^{(v)}}$. Thus $\cO_{\mathcal{X},\eta} = R_{\pi R}$, whose maximal ideal is generated by $\pi$. Consequently, $\tv(\pi) = 1$, i.e.\ $e(\tv|v) = 1$, and
\begin{equation}
L^{(\tv)} = R_{\pi R}/\pi R_{\pi R} = \operatorname{Frac}(R/\pi R),
\end{equation}
which is the function field $F^{(v)}(\mathcal{X}_{F^{(v)}})$ of the special fiber.

\vskip2mm

\noindent (iii) Since $\mathcal{X}_{F^{(v)}}$ is smooth of dimension $n$ over $F^{(v)}$, the sheaf of relative differentials $\Omega^1_{\mathcal{X}_{F^{(v)}}/F^{(v)}}$ is locally free of rank $n$ (see, e.g., \cite[Proposition 10.4]{Hartshorne}). Passing to the stalk at $\eta$ and using (ii), we see that $\Omega^1_{L^{(\tv)}/F^{(v)}}$ is an $L^{(\tv)}$-vector space of dimension $n = \operatorname{tr.deg}_{F^{(v)}}L^{(\tv)}$. This is equivalent to $L^{(\tv)}/F^{(v)}$ being separably generated (see \cite[Ch. II, \S17, Theorem 41]{Zariski-Samuel}).
\end{proof}

\vskip3mm

\noindent {\bf 2.5. Models of conics.} In this subsection, we describe explicit models of the function field of a conic, which will be used in \S\ref{S:SB} and \S\ref{S:FinUnram}. Let $K$ be a finitely generated field of characteristic $0$, fix $a, b \in K^\times$, and consider the projective conic $Q\colon ax^2 + by^2 = z^2$ in $\mathbb{P}^2_K$, as well as its affine open subset $Q_1 = Q \cap \{z \neq 0\}$ (i.e., $Q_1$ is the affine conic $ax^2 + by^2 = 1$). Set $L= K(Q) = K(Q_1)$ to be the corresponding function field.

Suppose $V$ is a divisorial set of places of $K$. In view of condition (*)
and Lemma \ref{L:AdeleRecall}, we may (and will) assume that $V$ arises from an affine model $\Spec A$ of $K$, where $A$ is a noetherian unique factorization domain, and that
$$
\operatorname{char} K^{(v)} \neq 2 \quad\text{and}\quad v(a) = v(b) = 0  \ \ \ \text{for all} \  v \in V.
$$

\begin{remark}\label{R:23}
Since $A$ is a noetherian unique factorization domain, in particular a normal noetherian domain, we have
$$
A = \bigcap_{{\rm ht}~\mathfrak{p} = 1} A_{\mathfrak{p}}
$$
by \cite[Theorem 11.5]{Matsumura}. Consequently, the condition $v(a) = v(b) = v(2) = 0$ for all $v \in V$ implies that $a, b$, and $2$ are units in $A$.
\end{remark}

Let us set \begin{equation}\label{eq:models}
\cQ := \operatorname{Proj}\big(A[x,y,z]/(ax^2+by^2-z^2)\big) \quad\text{and}\quad \cQ_1 := \Spec B, \quad B := A[x,y]/(ax^2+by^2-1),
\end{equation}
so that $\cQ_1 = \cQ\cap\{z \neq 0\}$. For $\mathfrak p \in \Spec A$ with residue field $\kappa(\mathfrak p) = A_{\mathfrak{p}}/\mathfrak{p} A_{\mathfrak{p}}$, we denote by $Q_{\mathfrak p}$ the conic $\bar ax^2 + \bar by^2 = z^2$ over $\kappa(\mathfrak p)$, where $\bar a, \bar b \in \kappa(\mathfrak p)^\times$ are the images of $a$ and $b$. If $v \in V$ corresponds to a prime $\mathfrak p$, we will write $Q^{(v)}= Q_{\mathfrak p}$. 

\begin{lemma}\label{lem23} \ \
\begin{enumerate}[label=\textup{(\roman*)}, leftmargin=*]

\vskip1mm

  \item[{\rm (i)}] $\cQ$ and $\cQ_1$ are normal integral schemes with function field $L$ that are smooth over $\Spec A$. Their fibers over $\mathfrak p \in \Spec A$ are $Q_{\mathfrak p}$ and $Q_{\mathfrak p}\cap\{z \neq 0\}$, respectively, which are smooth geometrically integral curves over $\kappa(\mathfrak p)$.
  \vskip1mm

  \item[{\rm (ii)}] Let $\tilde Q \subset Q$ be the open subset obtained by removing finitely many closed points $P_j$ \textup{(}$j \in J$\textup{)}, and denote by $Y_j$ the closure of $P_j$ in $\cQ$. Then, after possibly replacing $\Spec A$ by a suitable principal open subset, the $Y_j$ are pairwise disjoint, normal, and finite over $\Spec A$ \textup{(}so that $\bigcup_{j\in J} Y_j$ is a normal closed subscheme of $\cQ$ of pure codimension~$1$\textup{)}, and $\tilde\cQ := \cQ \setminus \bigcup_{j\in J} Y_j$ is a normal integral scheme that is smooth over $\Spec A$ and whose fibers are nonempty open subsets of the conics $Q_{\mathfrak p}$.
\end{enumerate}
\end{lemma}

\begin{proof} (i) The standard affine charts $\{z \neq 0\}$, $\{x \neq 0\}$, $\{y \neq 0\}$ of $\cQ$ have coordinate rings
\[
B \simeq A[x][y]/\big(y^2 - b^{-1}(1 - ax^2)\big), \qquad A[y][z]/\big(z^2 - (a + by^2)\big), \qquad A[x][z]/\big(z^2 - (ax^2 + b)\big)
\]
(for $B$, we use the fact that $b \in A^\times$). Each has the form $R[w]/(w^2 - g)$, with $R$ a polynomial ring in one variable over $A$ and $g \in R$. Hence, it is a free $R$-module with basis $\{1, w \}$, and therefore a free $A$-module. Consequently, $\cQ$ and $\cQ_1$ are flat and of finite presentation over $A$.

Next, since $a, b, 2 \in A^{\times}$, the Jacobian criterion implies that each $Q_{\mathfrak{p}}$ is a smooth conic over $\kappa(\mathfrak{p})$. Also, since a nondegenerate ternary quadratic form is irreducible over an algebraically closed field (see, e.g., \cite[Ch.~I, Exercise~5.12]{Hartshorne}), we see that $Q_{\mathfrak{p}}$ is geometrically integral. Moreover, $Q_{\mathfrak p}\cap\{z \neq 0\}$ is a nonempty open subset of $Q_{\mathfrak p}$ (since $Q_{\mathfrak p}\cap\{z = 0\}$ is clearly finite), hence is also smooth and geometrically integral.

Now, since $\cQ$ and $\cQ_1$ are flat and of finite presentation over $\Spec A$ and every fiber is smooth over the residue field, it follows that $\cQ$ and $\cQ_1$ are smooth over $\Spec A$ (see \cite[Tag 01V8]{Stacks}). Since $A$ is normal, applying \cite[Tag 033C]{Stacks} to the coordinate rings of the standard charts above, we conclude that $\cQ$ and $\cQ_1$ are normal.

Finally, let $C$ be the coordinate ring of one of the standard charts (e.g.\ $C = B$). Since $C$ is flat, hence torsion-free, over $A$, it embeds into $C\otimes_A K$, which is the coordinate ring of a nonempty affine open subset of the integral curve $Q$ over $K$. Thus, $C$ is an integral domain. Clearly, any two of the charts intersect, which implies that $\cQ$ is irreducible, and hence integral. The same is also true for its open subscheme $\cQ_1$. Both $\cQ$ and $\cQ_1$ obviously have function field $L$.

\vskip2mm

\noindent (ii)  Let $Y = \bigcup_{j\in J} Y_j$, with its reduced subscheme structure. Since $\operatorname{char} K = 0$, its generic fiber $\coprod_{j\in J}\Spec\kappa(P_j)$ is finite and \'etale over $K$.
Hence, after replacing $\Spec A$ by a suitable principal open subset, we may assume that $Y$ is finite and \'etale over $\Spec A$ (cf. \cite[Theorem 3.2.1]{Poonen}). Then $Y = \Spec C$ is affine (since it is finite over $\Spec A$ --- see \cite[Ch.~II, Exercise~3.4]{Hartshorne}), with
$C$ \'etale over $A$. Therefore, since $A$ is normal, $C$ is also normal by \cite[Tag 033C]{Stacks}.
Hence, the irreducible components $Y_j$ of $Y$ are normal and pairwise disjoint by \cite[Tag 0357]{Stacks} (see the proof of \cite[Th\'eor\`eme~6.2]{CTFinitude} for a similar argument). Finally, the fibers of $Y$ are finite, so the fibers of $\tilde\cQ$ are nonempty open subsets of the conics $Q_{\mathfrak p}$, and the remaining assertions follow from (i).

\end{proof}

Let us now consider the divisorial sets of places $\tilde V(\cQ)$ and $\tilde V(\cQ_1)$ of $L$ associated with the models $\cQ$ and $\cQ_1$, respectively. We recall that a prime divisor is called \emph{horizontal} if it dominates $\Spec A$ and is said to be \emph{vertical} otherwise (cf. \cite[Ch. 8, \S3.1]{Liu}).
We write
$$
\tilde V(\cQ) = \tilde V_0(\cQ) \cup \tilde V_1(\cQ) \ \ \ \ \text{and} \ \ \ \tilde V(\cQ_1) = \tilde V_0(\cQ_1) \cup \tilde V_1(\cQ_1)
$$
for the corresponding decompositions into \emph{horizontal} (or \emph{geometric}) and \emph{vertical} (or \emph{arithmetic}) places. We note that if $v \in V$ corresponds to a prime ideal $\mathfrak{p} \subset A$, so that $\cO_{K,v} = A_{\mathfrak p}$, then by Lemma~\ref{lem23}(i), the scheme $\cQ\times_A \cO_{K,v}$ over $\Spec\cO_{K,v}$ satisfies the assumptions of \S2.4, with generic fiber $Q$ and special fiber $Q^{(v)}$. Hence, the construction of \S2.4 yields a discrete valuation $\tv$ of $L$ extending $v$, with the valuation ring $\mathcal{O}_{L, \tilde{v}}$ being
the local ring of $\cQ$ at the generic point of $Q^{(v)}$.

\begin{lemma}\label{lem24} \ \
\begin{enumerate}[label=\textup{(\roman*)}, leftmargin=*]
 \vskip1mm

  \item \textup{(Horizontal places)} The horizontal prime divisors of $\cQ$ \textup{(}resp.\ $\cQ_1$\textup{)} are the closures in $\cQ$ \textup{(}resp.\ $\cQ_1$\textup{)} of the closed points of $Q$ \textup{(}resp.\ $Q_1$\textup{)}, and the corresponding places of $L$ restrict to zero on $K$.
   \item \textup{(Vertical places)} The vertical prime divisors of $\cQ$ are the closures in $\cQ$ of the fibers $Q^{(v)}$ for $v \in V$. Consequently, $\tilde V_1(\cQ)$ consists of the valuations $\tv$ \textup{(for }$v \in V$\textup{)} constructed above, and since each fiber $Q^{(v)}$ is irreducible, $\tv$ is the unique vertical place above $v$. Moreover, all vertical prime divisors of $\cQ$ meet $\cQ_1$, hence $\tilde V_1(\cQ_1) = \tilde V_1(\cQ)$.
  \item For $\tv \in \tilde V_1(\cQ)$ over $v \in V$, we have $e(\tv|v) = 1$ and $L^{(\tv)} = K^{(v)}(Q^{(v)})$.
  \item If $\tv \in \tilde V_0(\cQ_1)$ is a horizontal place that corresponds to the height one prime $\mathfrak P \subset B$, then $\cO_{L,\tv} = B_{\mathfrak P}$ contains $B$ and $K$. If $\tv \in \tilde V_1(\cQ_1)$ is a vertical place that lies above $v \in V$, with corresponding prime $\mathfrak p \subset A$, then $\mathfrak pB$ is a prime ideal of $B$ and $\cO_{L,\tv} = B_{\mathfrak pB}.$
 \end{enumerate}
\end{lemma}

\begin{proof}
(i), (ii) These are standard --- the reader can consult \cite[Ch.~8, Proposition~3.4]{Liu}, whose proof, based on the dimension formula for flat morphisms (cf. \cite[Theorem~15.1]{Matsumura}), applies without any changes in the present setting. In (ii), the uniqueness of the vertical place over $v$ comes from the irreducibility of $Q^{(v)}$ (Lemma~\ref{lem23}(i)), and the vertical prime divisors meet $\cQ_1$ because, as noted in the proof of Lemma~\ref{lem23}(i), $Q^{(v)}\cap\{z = 0\}$ is finite.

\vskip2mm

\noindent (iii) Since $\tv$ is the vertical valuation of $\cQ\times_A A_{\mathfrak p}$ over $v$ as discussed above, the statements follow immediately from
Lemma~\ref{lem22}.

\vskip2mm

\noindent (iv) In the horizontal case, $\mathfrak P$ lies over the generic point of $\Spec A$, so  $\mathfrak P \cap A = 0$. Hence,
$A\setminus\{0\} \subset B\setminus\mathfrak P$, and therefore $B_{\mathfrak P}$ contains $B$ and $K$. In the vertical case, recall from the proof of Lemma~\ref{lem23}(i) that $B$ is a free $A$-module. Hence $B/\mathfrak pB = B\otimes_A A/\mathfrak p$ is a free $A/\mathfrak p$-module, and therefore embeds into its localization $B\otimes_A\kappa(\mathfrak p)$ (note that $\kappa(\mathfrak p)$ is the fraction field of $A/\mathfrak p$). Since $B\otimes_A\kappa(\mathfrak p)$ is the coordinate ring of $Q^{(v)}\cap\{z\neq 0\}$, which is an integral domain by Lemma~\ref{lem23}(i), it follows that $\mathfrak pB$ is a prime ideal of $B$, namely the kernel of $B \to B\otimes_A\kappa(\mathfrak p)$; thus, it corresponds to the generic point of $Q^{(v)}$. Consequently, $\cO_{L,\tv} = B_{\mathfrak pB}$.
\end{proof}

\vskip3mm

\noindent {\bf 2.6. Unramified cohomology.}
Let $F$ be a field equipped with a discrete valuation $v$ and suppose $M$ is a finite $\Ga(F^{\rm sep}/F)$-module that is unramified at $v$, in the sense that for some (equivalently, any) extension $w$ of $v$ to $F^{\rm sep}$, the inertia subgroup
$I_w \subset \Ga(F^{\rm sep}/F)$ acts trivially on $M$.
Then there exists a residue map
$$
\rho^i_{M,v} \colon H^i(F, M) \to H^{i-1}(F^{(v)}, M(-1)) \ \ \ (i \geq 1)
$$
(cf. \cite[Chapter II]{GMS}). In particular, for every integer $\ell$ that is prime to ${\rm char}~F^{(v)}$, all integers $j$, and all $i \geq 1$, we have residue maps
$$
\rho^i_{\ell, j, v} \colon H^i(F, \mu_{\ell}^{\otimes j}) \to H^{i-1}(F^{(v)}, \mu_{\ell}^{\otimes (j-1)}),
$$
where $\mu_{\ell}$ is the group of $\ell$-th roots of unity and $\mu_{\ell}^{\otimes j}$ denotes the $j$-th Tate twist of $\mu_{\ell}$. In this paper, we will only need to consider the case where $\ell$ is a power of 2, and then, to avoid cumbersome notation, we will denote the residue maps that arise simply by $\rho^i_v$ whenever this does not lead to confusion.

We say that an element of $H^i(F,M)$ is {\it unramified at $v$} if it lies in the kernel of the relevant residue map. Furthermore, given a set $V$ of discrete valuations of $F$ such that the residue maps exist for all $v \in V$, one defines, for all $i \geq 1$, the corresponding unramified cohomology groups by
$$
H^i(F, M)_V = \bigcap_{v \in V} \ker \rho^i_{M,v}.
$$

We will need to consider one particular extension of the preceding constructions to the case of infinite torsion coefficients. Namely, suppose $p$ is a prime $\neq {\rm char}~F$. For any integer $j$, we set
$$
\Q_p/\Z_p (j) = \lim_{\longrightarrow} \mu_{p^k}^{\otimes j},
$$
where the limit is taken with respect to the natural maps
\begin{equation}\label{E:IndLimit}
\mu_{p^m}^{\otimes j} \to \mu_{p^n}^{\otimes j} \ \ \ \text{for} \ \  n \geq m.
\end{equation}
Then for all $i \geq 1$, we define
$$
H^i(F, \Q_p/\Z_p(j)) = \lim_{\longrightarrow} H^i(F, \mu_{p^k}^{\otimes j}),
$$
with the limit taken with respect to the maps induced by (\ref{E:IndLimit}). Now, if $V$ is a set of discrete valuations such that ${\rm char}~F^{(v)} \neq p$ for all $v \in V$, then the residue maps
$$
\rho^i_{p^k, j, v} \colon H^i(F, \mu_{p^k}^{\otimes j}) \to H^{i-1}(F^{(v)}, \mu_{p^k}^{\otimes (j-1)})
$$
assemble to yield residue maps
$$
\rho^i_{p^{\infty}, j, v} \colon H^i(F, \Q_p/\Z_p(j)) \to H^{i-1} (F^{(v)}, \Q_p/\Z_p(j-1))
$$
for all $v \in V$ and $i \geq 1.$ As above, we then can consider the subgroup of unramified elements
$$
H^i(F, \Q_p/\Z_p(j))_V := \bigcap_{v \in V} \ker \rho^i_{p^{\infty}, j, v}.
$$

We would also like to mention a localization result for unramified cohomology that will be used several times later on. Suppose that $X$ is an excellent integral noetherian normal scheme with function field $K$ and $Y \subset X$ is a {\it normal} closed subscheme of pure codimension 1. Thus, writing the decomposition $Y = \bigcup_{j \in J} Y_j$ of $Y$ into its irreducible components, each $Y_j$ has codimension 1 in $X$, and we let $\kappa_j$ denote the function field of $Y_j$. Set $U = X \setminus Y$, and let $n$ be a positive integer invertible on $X$. Then for any $i \geq 1$ and any $\ell$, we  have a canonical embedding
$$
s \colon H^i(K , \mu_n^{\otimes \ell})_{V_X} \longrightarrow H^i(K , \mu_n^{\otimes \ell})_{V_U},
$$
where $V_X$ (resp. $V_U$) is the set of discrete valuations of $K$ associated with the prime divisors of $X$ (resp. of $U$).
On the other hand, for each $j \in J$, we have the residue map $H^i(K , \mu_n^{\otimes \ell}) \to H^{i-1}(\kappa_j , \mu_n^{\otimes (\ell-1)})$, and we let
$$
r \colon H^i(K , \mu_n^{\otimes \ell}) \longrightarrow \bigoplus_{j \in J} H^{i-1}(\kappa_j , \mu_n^{\otimes (\ell-1)})
$$
denote the direct sum of all these maps. We then have the following.
\begin{prop}\label{P:Open} {\rm (see \cite[Proposition 6.4]{Spinor})}
The above maps $r$ and $s$ give the following exact sequence of unramified cohomology groups
\begin{equation}\label{E:Open}
0 \to H^i(K , \mu_n^{\otimes \ell})_{V_X} \stackrel{s}{\longrightarrow} H^i(K , \mu_n^{\otimes \ell})_{V_U} \stackrel{r}{\longrightarrow} \bigoplus_{j \in J} H^{i-1}(\kappa_j , \mu_n^{\otimes (\ell-1)})_{V_{Y_j}}.
\end{equation}
\end{prop}



\section{Morita equivalence for skew-hermitian forms over quaternion algebras}\label{S:Morita}

Suppose $F$ is a field of characteristic $\neq 2$ and $D$ is a quaternion algebra over $F$. If $D \simeq M_2(F)$ is the split algebra equipped with the standard involution
$$
X \mapsto J X^t J^{-1}, \ \ \ \text{where} \ \ J = \left(\begin{array}{cr} 0 & -1 \\ 1 & 0 \end{array} \right)
$$
then by Morita equivalence, one can associate to a nondegenerate $n$-dimensional skew-hermitian form $h$ over $D$ a nondegenerate $2n$-dimensional quadratic form $q_h$ over $F$ (see \cite[Chapter 1, \S9]{Knus} for the general theory of Morita equivalence and \cite[Ch. 10, \S3]{Scharlau} for a discussion of skew-hermitian forms over quaternion algebras). By general properties of Morita equivalence, it is known that $h$ is completely determined by $q_h$, in the sense that two skew-hermitian forms are equivalent if and only if the corresponding quadratic forms are equivalent.
In Proposition~\ref{prop32I} below, we give an explicit description of this correspondence that will be needed for our purposes.

We begin by fixing our setup. Let $D$ be a quaternion division algebra over $F$ with its canonical involution
$$
\iota \colon D \to D, \ \ \ z \mapsto \bar{z}.
$$
Following the conventions of \cite{Knus} and \cite{Scharlau}, we will work with {\it right} modules over $D$ (as well as over scalar extensions of $D$). Thus, if $M$ is a right $D$-module, then a skew-hermitian form
$
h \colon M \times M \to D
$
is a bi-additive map that satisfies
$$
h (x\lambda, y\mu) = \bar{\lambda} h(x,y)\mu \ \ \ \text{and}  \ \ \ h(y,x) = -\overline{h(x,y)} \ \ \ \text{for all} \ x,y \in M \ \text{and} \ \lambda, \mu \in D.
$$
If $h$ is a skew-hermitian form on $D^n$ and $\{ e_1, \dots, e_n \}$ is a basis of $D^n$, we consider the Gram matrix
$$
C = (h(e_l, e_m)) \in M_n(D)
$$
of $h$ with respect to this basis. Writing $C^* = (\bar{c}_{ml}) = (\bar{c}_{lm})^t$, we have $C^* = -C$ and $h(x,y) = x^* C y$ (in coordinates with respect to $\{ e_1, \dots, e_n \}$), and, moreover, if $\{f_1, \dots, f_n \}$ is another basis with
$$
f_m = \sum_{l = 1}^n e_l p_{lm},
$$
then the corresponding Gram matrix is $P^* C P$, where $P = (p_{lm}).$

Next, if $F'/F$ is a field extension, then $\iota$ extends to an involution $\iota \otimes {\rm id}$ on $D \otimes_F F'$ through its action on the first factor. By a slight abuse of notation, we will still write $z \mapsto \bar{z}$ for $z \in D \otimes F'$. The base change $h_{F'} = h \otimes_F F'$ is a skew-hermitian form on $D^n \otimes_F F' = (D \otimes_F F')^n$ and $C$ is the Gram matrix of $h_{F'}$ with respect to the basis $\{ e_1 \otimes 1, \dots, e_n \otimes 1 \}.$

Suppose now that $F'$ is a splitting field of $D$ and fix an isomorphism $\rho \colon D \otimes_F F' \to M_2(F')$ of $F'$-algebras. This gives rise to an $F'$-algebra isomorphism
$$
\Phi_{\rho} \colon M_n(D \otimes_F F') \to M_{2n}(F') \ \ \ \text{with} \ \ \ \Phi_{\rho}(X) = (\rho(x_{lm}))_{l, m} \ \ \ \text{for} \ \ X = (x_{lm}) \in M_n(D \otimes F').
$$
We will also need the block diagonal matrix $\mathbf{J} := \operatorname{diag}(J,\dots,J) \in M_{2n}(F')$.

\begin{prop}\label{prop32I}
With notations as above, we have the following:
\begin{enumerate}[label=\textup{(\roman*)}, leftmargin=*, itemsep=0.3em]
  \item For all $z \in D\otimes_F F'$, we have
  \[
  \rho(\bar z) = J\rho(z)^tJ^{-1}.
  \]
  \item If $H$ is a nondegenerate skew-hermitian form on $(D\otimes_F F')^n$ with Gram matrix $C$, then the matrix
  \[
  S_\rho(C) := \mathbf{J}^{-1}\Phi_\rho(C) \in M_{2n}(F')
  \]
  is symmetric and invertible, and hence defines a nondegenerate quadratic form $q_{H,\rho}$ on $(F')^{2n}$.
  In particular, if $h$ is a nondegenerate skew-hermitian form on $D^n$ and $C$ is a Gram matrix of $h_{F'}$ \textup{(}for instance, a Gram matrix of $h$\textup{)}, then we obtain in this way a nondegenerate quadratic form $q_{h_{F'},\rho}$ on $(F')^{2n}$.
  \item If $P \in \GL_n(D\otimes_F F')$, so that $P^*CP$ is the Gram matrix of $H$ with respect to another basis, then
  \[
  S_\rho(P^*CP) = \Phi_\rho(P)^t\,S_\rho(C)\,\Phi_\rho(P).
  \]
  Consequently, up to equivalence, the quadratic form $q_{H,\rho}$ does not depend on the choice of the Gram matrix $C$ of $H$.
  \item Two nondegenerate skew-hermitian forms $H$ and $H'$ on $(D\otimes_F F')^n$ are equivalent if and only if the corresponding quadratic forms $q_{H,\rho}$ and $q_{H',\rho}$ are equivalent.
\end{enumerate}
\end{prop}

\begin{proof} \ \

\noindent (i) We recall that the canonical involution of $D \otimes_F F'$ is the unique $F'$-linear involution satisfying $z \bar{z} \in F'$ for all $z \in D \otimes_F F'$. Since $\rho$ preserves this property, it carries the canonical involution of $D\otimes_F F'$ to the standard involution of $M_2(F')$, which is the adjugate (see \cite[Chapter 3]{Vo} for a detailed discussion of involutions).

\vskip3mm

\noindent (ii) We first observe that using (i) and direct computation, we have
\begin{equation}\label{E:Morita-matrix-formula}
\Phi_\rho(X^*) = \mathbf{J}\,\Phi_\rho(X)^t\,\mathbf{J}^{-1} \ \ \ \ \text{for all} \ \ X \in M_n(D\otimes_F F').
\end{equation}
We also note that $\mathbf{J}^t = \mathbf{J}^{-1} = -\mathbf{J}$. Consequently, since $\Phi_\rho(C)^t = \mathbf{J}^{-1}\Phi_\rho(C^*)\mathbf{J}$ by (\ref{E:Morita-matrix-formula}) and $C^* = -C$ by assumption, we obtain
$$
S_\rho(C)^t = \Phi_\rho(C)^t\,\mathbf{J} = \mathbf{J}^{-1}\Phi_\rho(C^*)\,\mathbf{J}^2 = -\mathbf{J}^{-1}\Phi_\rho(C^*) = \mathbf{J}^{-1}\Phi_\rho(C) = S_\rho(C),
$$
showing that $S_{\rho}(C)$ is symmetric. Furthermore, since $H$ is nondegenerate, so that $C$ is invertible, the matrix $\Phi_\rho(C)$ is invertible; hence, since $\det\mathbf{J} = 1$, it follows that $S_{\rho}(C)$ is invertible.

\vskip3mm

\noindent (iii) Let $P \in \GL_n(D\otimes_F F')$. Then using the multiplicativity of $\Phi_\rho$, together with (\ref{E:Morita-matrix-formula}), we obtain
$$
S_\rho(P^*CP) = \mathbf{J}^{-1}\,\mathbf{J}\,\Phi_\rho(P)^t\,\mathbf{J}^{-1}\,\Phi_\rho(C)\,\Phi_\rho(P) = \Phi_\rho(P)^t\,S_\rho(C)\,\Phi_\rho(P).
$$

\vskip3mm

\noindent (iv) Let $C$ and $C'$ be Gram matrices of $H$ and $H'$. If $H$ and $H'$ are equivalent, then $C' = P^*CP$ for some $P \in \GL_n(D\otimes_F F')$, and hence $q_{H,\rho}$ and $q_{H',\rho}$ are equivalent by (iii). Conversely, suppose $q_{H,\rho}$ and $q_{H',\rho}$ are equivalent, so that
$$
S_\rho(C') = N^tS_\rho(C)N \ \ \  \text{for some} \  N \in \GL_{2n}(F').
$$
Set $P = \Phi_\rho^{-1}(N) \in \GL_n(D\otimes_F F')$. Then by (iii), we have $S_\rho(P^*CP) = S_\rho(C')$, and hence $P^*CP = C'$. Thus, $H$ and $H'$ are equivalent.

\end{proof}

\begin{remark}\label{rem32I}
We note that in the context of the general theory described in \cite[Chapter 1, \S9]{Knus} and \cite[Ch. 10, \S3]{Scharlau},
$q_{H,\rho}$ is precisely the quadratic form attached to $H$ by Morita equivalence via $\rho$, with respect to the alternating form $b_0(u,w) = u^tJ^{-1}w$ on $(F')^2$, whose adjoint involution is $X \mapsto JX^tJ^{-1}$.
\end{remark}

In the proof of Theorem \ref{T-SB}, we will apply Proposition \ref{prop32I} in the context of the generic splitting of a quaternion division algebra by the function field of its Severi-Brauer conic. More precisely, let $D$ be a quaternion division algebra over $F$ with a fixed presentation $D = \displaystyle{\left( \frac{a,b}{F} \right)}.$ Thus, $D$ has a standard $F$-basis $1, i, j, k$ satisfying
$$
i^2 = a, \ \ j^2 = b, \ \ ij =k, \ \ ij = -ji.
$$
Also, denote by $F(Q)$ the function field of the corresponding Severi-Brauer conic $Q$ --- in other words, $F(Q)$ is simply the fraction field
of the integral domain $F[x,y]/(ax^2 + by^2 -1)$. We have the following elementary statement, which is easily verified by direct calculation.



\begin{lemma}\label{L:SplittingSB}
Let $D = \displaystyle{\left(\frac{a,b}{F} \right)}$ be a quaternion division algebra over a field $F$ of characteristic $\neq 2$, and denote by $F(Q)$ the function field of the corresponding affine Severi-Brauer conic $ax^2 + by^2 = 1.$ Then we have an isomorphism
$$
\varphi \colon D \otimes_F F(Q) \stackrel{\sim}{\longrightarrow} M_2(F(Q))
$$
given by
$$
\varphi(i) = \left( \begin{array}{cc} ax & -y \\ -aby & -ax \end{array} \right) \ \ \ \text{and} \ \ \ \varphi(j) = \left( \begin{array}{rr} by & x \\ abx & -by \end{array} \right).
$$
\end{lemma}

\vskip2mm

In the rest of the paper, given a skew-hermitian form $h$ over $D$, we will write $q_{h_{F(Q)}} = q_{h_{F(Q)},\varphi}$ for the quadratic form over $F(Q)$ corresponding to $h_{F(Q)} = h\otimes_F F(Q)$ via the splitting $\varphi$ of Lemma~\ref{L:SplittingSB}. By Proposition~\ref{prop32I}(ii), if $h$ has Gram matrix $C$, then $q_{h_{F(Q)}}$ is the quadratic form with Gram matrix $\mathbf{J}^{-1}\Phi_\varphi(C)$.

\vskip3mm

To conclude this section, we would like to mention the following result, which will play an important role in the proof of Theorem \ref{T-SB}.

\begin{prop}\label{P-PSS} {\rm (\cite[Proposition 3.3]{PSS})} Let $D$ be a quaternion algebra over a field $F$ and let $\iota$ be the canonical involution on $D$. Let $Q$ be the Severi-Brauer conic associated to $D$ and let $F(Q)$ be the function field of $Q$. Then the canonical homomorphism
$$
W^{-1}(D, \iota) \to W^{-1}(D \otimes_F F(Q), \iota \otimes 1)
$$
is injective, where $W^{-1}$ is the Witt group of skew-hermitian forms.
\end{prop}

\section{Quadratic and skew-hermitian forms with good reduction}\label{S:quadratic-hermitian-forms-good-reduction}

The purpose of this section is to develop some notions
pertaining to quadratic and skew-hermitian forms with good reduction and also to establish a conditional finiteness result (Proposition \ref{P:FormsGoodReduction}) that will be needed for
the proof of Theorem \ref{T-SB}.

Throughout this section, we will work over a field $F$ equipped with a discrete valuation $v$, and we will assume that the residue field $F^{(v)}$ has characteristic $\neq 2$ (so that 2 is a unit in the valuation ring $\mathcal{O}_{F_v} \subset F_v$ of the completion) --- in particular, ${\rm char}~F \neq 2$ as well. Given a free $\mathcal{O}_{F_v}$-module $M$ of finite rank, we recall that a symmetric bilinear form $b \colon M \times M \to \mathcal{O}_{F_v}$ is {\it regular} if the canonical induced map $M \to M^{\vee} = \Hom_{\mathcal{O}_{F_v}}(M, \mathcal{O}_{F_v})$ is an isomorphism. Furthermore, a quadratic form $q \colon M \to \mathcal{O}_{F_v}$ is said to be regular if its associated symmetric bilinear form is regular. If ${\rm rank}~M = n$, we will say that $q$ has rank $n$.
As usual, we let $\langle a \rangle \colon \mathcal{O}_{F_v} \to \mathcal{O}_{F_v}$ denote the standard quadratic form $x \mapsto ax^2.$ Then, in view of our assumption that $2 \in \mathcal{O}_{F_v}^{\times}$, we have the following basic fact.

\begin{lemma}\label{L:RegularQForm}
Any regular quadratic form $(M,q)$ over $\mathcal{O}_{F_v}$ has an orthonormal basis
$$
(M, q) \simeq \langle a_1\rangle \perp \dots \perp \langle a_n \rangle
$$
for $a_i \in \mathcal{O}_{F_v}^{\times}$, where $n = {\rm rank}~M.$ In other words, $q$ is equivalent over $\mathcal{O}_{F_v}$ to a diagonal form $a_1 x_1^2 + \cdots + a_n x_n^2.$
\end{lemma}
\begin{proof}
This is a special case of \cite[Proposition 3.4]{Baeza}.
\end{proof}

We now make the following definition.

\vskip2mm

\begin{definition}\label{D:GRquadratic}
A nondegenerate $n$-dimensional quadratic form $q$ over $F$ has {\it good reduction at $v$} (or is {\it unramified at $v$}) if there exists a regular quadratic form $\tilde{q}$ of rank $n$ over $\mathcal{O}_{F_v}$ such that
$$
\tilde{q} \otimes_{\mathcal{O}_{F_v}} F_v \simeq q \otimes_F F_v.
$$
\end{definition}

\vskip2mm

With this terminology, we have the following statement. 

\begin{lemma}\label{L-SpinGR} {\rm (cf. \cite[Example 2.3]{CRR-MaxTori})}
Suppose $q$ is an $n$-dimensional nondegenerate quadratic form over $F$ and let $G = {\rm Spin}_n(q)$ be the corresponding spinor group. Denote by $q_v = q \otimes_F F_v$ the base change of $q$ to $F_v$. Then $G$ has good reduction at $v$ if and only if $\lambda_v q_v$ has good reduction at $v$ for some $\lambda_v \in F_v^{\times}.$

\end{lemma}

Next, let us discuss the analogous notions for skew-hermitian forms over quaternion algebras. Let $D$ be a quaternion division algebra over $F$ and assume that $D$ is unramified at $v$, i.e. the class $[D] \in {}_2\Br(F) \simeq H^2(F, \mu_2)$ is contained in the kernel of the residue map
$$
\rho^2_v \colon H^2(F, \mu_2) \to H^1(F^{(v)}, \mu_2).
$$
Set $D_v = D \otimes_F F_v.$ Then either $D_v$ is isomorphic to $M_2(F_v)$ or $D_v$ is a quaternion division algebra over $F_v$. In the first case,
we {\it fix} an isomorphism $\psi\colon D_v \stackrel{\sim}{\longrightarrow} M_2(F_v)$ and set $\mathcal{A}_v = \psi^{-1}(M_2(\cO_{F_v}))$ (which is clearly an Azumaya algebra over $\cO_{F_v}$).

In the second case, it is well-known that the valuation $v$ on $F_v$ extends uniquely to a valuation $v_{D_v}$ on $D_v$ (see \cite[Corollary 2.2]{W2}), and we set
$$
\mathcal{A}_v = \{a \in D_v^{\times} \mid v_{D_v}(a) \geq 0 \} \cup \{ 0 \}
$$
to be the valuation ring of $v_{D_v}.$ Since $D$ is unramified at $v$, $\mathcal{A}_v$ is an Azumaya algebra over $\mathcal{O}_{F_v}$ with the property that $D_v \simeq \mathcal{A}_v \otimes_{\mathcal{O}_{F_v}} F_v$ (see \cite[Theorem 10.3]{Saltman} and \cite[Theorem 3.2]{W2}). Thus, in both
cases,
we have an Azumaya $\mathcal{O}_{F_v}$-algebra $\mathcal{A}_v \subset D_v$ such that
$$
D_v \simeq \mathcal{A}_v \otimes_{\mathcal{O}_{F_v}} F_v,
$$
and furthermore,
the canonical involution on $D_v$ induces an involution of $\mathcal{A}_v$.

Next, suppose $M$ is a free right $\mathcal{A}_v$-module of finite rank. We recall that a skew-hermitian form
$$
h \colon M \times M \to \mathcal{A}_v
$$
is called {\it regular} (or {\it nonsingular}) if the canonical induced map
$$
M \to M^{\vee} = \Hom_{\mathcal{A}_v}(M, \mathcal{A}_v), \ \ \ x \mapsto h(x, \cdot)
$$
is an isomorphism.

\vskip2mm

\begin{definition}\label{D:GRskew}
Suppose $h$ is an $n$-dimensional nondegenerate skew-hermitian form over a quaternion division algebra $D$ over $F$,
and denote by $h_v := h\otimes_F F_v$ the skew-hermitian form on $D_v^n$ obtained by extension of scalars. We say that $h$ has \emph{good reduction at $v$} if
\begin{enumerate}[label=\textup{(\roman*)}, leftmargin=*]
  \item $D$ is unramified at $v$; and
  \item there exist a free right $\mathcal{A}_v$-module $M$ of rank $n$ and a regular skew-hermitian form
  $$
  \tilde h\colon M\times M \to \mathcal{A}_v
  $$
  such that there is an equivalence
  $$
  \theta \colon (M \otimes_{\mathcal{O}_{F_v}} F_v, \tilde{h}_{F_v}) \stackrel{\sim}{\longrightarrow} (D_v^n, h_v)
  $$
  of skew-hermitian forms over $D_v$, where $\tilde{h}_{F_v}$ is the extension of $\tilde{h}$ to the right $D_v$-module $M\otimes_{\mathcal{O}_{F_v}}F_v$.

\end{enumerate}
\end{definition}

Note that condition (ii) of Definition \ref{D:GRskew} simply means that $h_v$ has a Gram matrix in $\mathrm{GL}_n(\mathcal{A}_v)$, namely the Gram matrix $C = (\tilde{h}(f_l, f_m))$ of $\tilde h$ with respect to an $\mathcal{A}_v$-basis $f_1,\dots,f_n$ of $M$, which is also the Gram matrix of $h_v$ with respect to the basis $\theta(f_1\otimes1),\dots,\theta(f_n\otimes1)$.

\vskip2mm

\begin{remark}\label{rem45I}
When $D_v$ is the split algebra, then the order $\mathcal{A}_v$, and hence the notion of good reduction in Definition \ref{D:GRskew}, depends on the choice of the isomorphism $\psi\colon D_v \stackrel{\sim}{\longrightarrow} M_2(F_v)$. However, we note that this dependence only involves a scaling factor.
Indeed, by the Skolem-Noether theorem, any other such isomorphism has the form $\psi\circ{\rm Int}(g)^{-1}$, with $g \in D_v^\times$, and it replaces $\mathcal{A}_v$ by $g\mathcal{A}_vg^{-1}$. If $C = (c_{lm}) \in \GL_n(\mathcal{A}_v)$ is a Gram matrix of $h_v$ (see the comment following Definition \ref{D:GRskew}),
then changing the basis by $P = g^{-1}I_n$ yields the Gram matrix $P^*CP$ of $h_v$. Since $\bar g^{-1} = \Nrd_{D_v/F_v}(g)^{-1}g$, the matrix $(gc_{lm}g^{-1}) = \Nrd_{D_v/F_v}(g)\,P^*CP$ is a Gram matrix of $\Nrd_{D_v/F_v}(g)\,h_v$ whose entries
lie in $g\mathcal{A}_vg^{-1}$. It follows that
$(h_v, D_v)$ has good reduction with respect to the order $\mathcal{A}_v$ if and only if $(\Nrd_{D_v/F_v}(g)\,h_v, D_v)$ has good reduction with respect to the order $g\mathcal{A}_vg^{-1}$.
\end{remark}



\vskip2mm

Using considerations similar to those employed in the proof of Lemma \ref{L-SpinGR} (see also \cite[Remark 4.3]{Srimathy}), one shows the following.


\begin{lemma}\label{L-SUGR}
Let $F$ be a field equipped with a discrete valuation $v$ such that ${\rm char}~F^{(v)} \neq 2$ and $D$ be a finite-dimensional quaternion division algebra over $F$.
Suppose $h$ is a finite-dimensional nondegenerate skew-hermitian form over
$D$, let $G = {\rm SU}(h, D)$ be the corresponding special unitary group, and set $\tilde{G}$ to be its universal cover. Then $\tilde{G}$ has good reduction at $v$ if and only if the form $\lambda_v h_v$ over $D_v$ has good reduction at $v$ for some $\lambda_v \in F_v^{\times}.$
\end{lemma}

In view of Remark~\ref{rem45I}, we note that the condition in Lemma~\ref{L-SUGR}, which asks for \emph{some} $\lambda_v \in F_v^\times$, does not depend on the choice of the isomorphism $\psi \colon D_v \stackrel{\sim}{\longrightarrow} M_2(F_v)$ in the split case.

\vskip2mm

In the remainder of this section, our objective will be to establish the next result, which relates the finiteness properties of quadratic forms having good reduction with respect to a set of discrete valuations of a field to those of the corresponding unramified cohomology groups. Although the argument that we give is an adaptation of the one used in the proof of \cite[Theorem 2.1]{Spinor}, since the required statement did not appear explicitly in {\it loc. cit.}, we provide full details for the sake of completeness.

\begin{prop}\label{P:FormsGoodReduction}
Suppose $F$ is a field of characteristic $\neq 2$ that is equipped with a set $V$ of discrete valuations such that ${\rm char}~F^{(v)} \neq 2$ for all $v \in V$. Let $n \geq 2$ be an integer and assume that the unramified cohomology groups $H^i(F, \mu_2)_V$ are finite for all $i = 1, \dots, \ell : = [\log_2 n]+ 1.$ Then the number of equivalence classes of nondegenerate $n$-dimensional  quadratic forms over $F$ having good reduction at all $v \in V$ is bounded above by $\displaystyle{
\prod_{i = 1}^{\ell} \vert H^i(F, \mu_2)_V \vert}
$ (in particular, it is finite).
\end{prop}

Before giving the proof, we recall several facts about the Witt ring and the Milnor isomorphism (see \cite[\S3]{Spinor} for more details). As above, we take $F$ to be a field of characteristic $\neq 2$. Let $W(F)$ and $I(F)$ denote the Witt ring of $F$ and its fundamental ideal, respectively. For a nondegenerate quadratic form $q$ over $F$, we denote by $[q]$ the corresponding class in $W(F)$. We recall that, as a consequence of Voevodsky's proof of the Milnor conjecture, there are natural isomorphisms of abelian groups
$$
\gamma_{F,d} \colon I(F)^d/I(F)^{d+1} \to H^d(F, \mu_2) \ \ \ \text{for all} \ d \geq 1
$$
(cf.\cite{OVV}). Now suppose that $F$ is complete with respect to a discrete valuation $v$, and let
$$
U(F) = \{u \in F \mid v(u) = 0 \}
$$
be the group of units of $F$. We set $W_0(F)$ to be the subring of $W(F)$ generated by the classes of $\langle u \rangle$ for $u \in U(F)$. We have the following statement.

\begin{lemma}\label{L:unram111} {\rm (cf. \cite[Lemma 3.3]{Spinor})}
Let $F$ be a field complete
with respect to a discrete valuation $v$ such that $\mathrm{char}\: F^{(v)} \neq 2$, and let $d \geq 1$. If $q$ is a nondegenerate quadratic form over $F$ such that $[\lambda q] \in I(F)^d \cap W_0(F)$ for some $\lambda \in F^{\times}$, then $[q] \in I(F)^d$ and the cohomology class
$\gamma_{F , d}([q]) \in H^d(F , \mu_2)$ is unramified at $v$.
\end{lemma}

\vskip2mm

\noindent {\it Proof of Proposition \ref{P:FormsGoodReduction}.}
Let $\{q_i\}_{i \in \cJ}$ be a family of $n$-dimensional nondegenerate quadratic forms over $F$ that have good reduction at all $v \in V$ and such that for $i, j \in \cJ$, $i \neq j$, the forms $q_i$ and $q_j$ are not equivalent. We will show that $\cJ$ is finite with
$$
\vert \cJ \vert \leq d_1 \cdots d_{\ell},
$$
where $d_i = \vert H^i(F, \mu_2)_V \vert$ for $i = 1, \dots, \ell.$

Fix $j_0 \in \cJ$. Then for any $j \in \cJ$ and any $v \in V$, the class of the form
$$
q_{j} \bot (-q_{j_0})
$$
lies in $W_0(F_v) \cap I(F_v)$ (see Lemma \ref{L:RegularQForm}). Hence, using Lemma \ref{L:unram111} with $d=1$, we see that
$$
\gamma_{F_v, 1}([q_j] - [q_{j_0}]) \in H^1(F_v, \mu_2)
$$
is unramified at $v$. Since this is true for all $v \in V$, we conclude that
$$
\gamma_{F,1}([q_j]-[q_{j_0}]) \in H^1(F, \mu_2)_V.
$$
Then one can find a subset $\cJ_1 \subset \cJ$ of size
$$
\geq \vert \cJ \vert/d_1
$$
(if $\cJ$ is infinite, this simply means that $\cJ_1$ is also infinite) such that the elements $[q_j] - [q_{j_0}] \in I(F)$ for $j \in \cJ_1$ have the same image under $\gamma_{F,1}$. Fix $j_1 \in \cJ_1$. Then for any $j \in \cJ_1$, we have
$$
\gamma_{F,1}([q_j] - [q_{j_1}]) = \gamma_{F,1}(([q_j] - [q_{j_0}]) - ([q_{j_1}] - [q_{j_0}])) = 0,
$$
which implies that $[q_j] - [q_{j_1}] \in I(F)^2.$ Consequently, for each $v \in V$, the class of the form $q_j \bot (-q_{j_1})$ lies in $W_0(F_v) \cap I(F_v)^2$, so Lemma \ref{L:unram111} with $d=2$ yields
$$
\gamma_{F,2}([q_j] - [q_{j_1}]) \in H^2(F, \mu_2)_V \ \ \ \text{for all} \ j \in \cJ_1.
$$
Then there exists a subset $\cJ_2 \subset \cJ_1$ of size $\geq \vert \cJ \vert/(d_1 d_2)$ such that the elements $[q_j] - [q_{j_1}]$ for $j \in \cJ_2$ have the same image under $\gamma_{F, 2}$. Consequently, fixing $j_2 \in \cJ_2$ and arguing as above, we obtain $[q_j]-[q_{j_2}] \in I(F)^3$ for all $j \in \cJ_2$. Inductively, we construct a nested sequence of subsets
$$
\cJ \supset \cJ_1 \supset \cJ_2 \supset \cdots \supset \cJ_{\ell}
$$
such that

\vskip1mm

\noindent (i) $\vert \cJ_m \vert \geq \vert \cJ \vert / (d_1 \cdots d_m)$; and
\vskip1mm

\noindent (ii) fixing $j_m \in \cJ_m$, we have $[q_j] - [q_{j_m}] \in I(F)^{m+1}$ for all $j \in \cJ_m$,

\vskip1mm

\noindent for any $m = 1, \dots, \ell.$ However, by \cite[Ch. X, Hauptsatz 5.1]{Lam}, the dimension of any positive-dimensional anisotropic form in $I(F)^{\ell + 1}$ is $\geq 2^{\ell+1} > 2^{\log_2 n + 1} = 2n$. So, the fact that
$$
[{q}_j] - [{q}_{j_{\ell}}] = [{q}_j \bot (- {q}_{j_{\ell}})] \in I(F)^{\ell + 1}
$$
implies that the form ${q}_j \bot (- {q}_{j_{\ell}})$ is hyperbolic, and consequently ${q}_j$ and ${q}_{j_{\ell}}$ are equivalent. Since $j \in \cJ_{\ell}$ was arbitrary and the forms ${q}_j$ for $j \in \cJ_{\ell}$ are assumed to be
pairwise inequivalent, we see that $\cJ_{\ell}$ actually reduces to a single element. Then the inequality in (i) yields the required estimation
$$
\vert \cJ \vert \leq d_1 \cdots d_{\ell},
$$
which concludes the proof.  \hfill $\Box$

\section{Proofs of Theorems \ref{T-SB} and \ref{T-GenusSpinor}: Reduction to Theorem \ref{T-SBUnram}}\label{S:SB}

Our goal in this section and the next one is to prove Theorem \ref{T-SB}. In this section, we reduce Theorem \ref{T-SB} to Theorem \ref{T-SBUnram} for a particular choice of divisorial set. We also deduce Theorem \ref{T-GenusSpinor} from Theorem \ref{T-SB}. The proof of Theorem \ref{T-SBUnram} is then given in \S\ref{S:FinUnram}, thereby completing the argument.

To fix notations, we take $K = k(C)$ to be the function field of a smooth geometrically integral curve $C$ over a number field $k$ and let $V$ be a divisorial set of places of $K$. Since $V$ obviously satisfies condition (*) introduced in \S2.2, without loss of generality, we can (and will) assume that the residue fields $K^{(v)}$ have characteristic $\neq 2$ for all $v \in V$.

Recall that our goal in Theorem \ref{T-SB} is to show that the set of $K$-isomorphism classes of the universal covers $\tilde{G}$ of special unitary groups $G = \mathrm{SU}_n(h,D)$, where $h$ is a nondegenerate $n$-dimensional skew-hermitian form over some quaternion division $K$-algebra $D$ with canonical involution,
having good reduction at all $v \in V$ is finite.
Now, by Lemma \ref{L-SUGR}, if $\tilde{G}$ has good reduction at all $v \in V$, then the quaternion algebra $D$ is unramified at all $v \in V$. On the other hand, since the 2-torsion of the unramified Brauer group
$$
{}_2\Br(K)_V := \bigcap_{v \in V} \ker \left( \rho_v^2 \colon H^2(K, \mu_2) \to H^1(K^{(v)}, \mu_2) \right)
$$
is finite by \cite[Theorem 1]{CRR3} (where $\rho_v^2$ are the residue maps for $v \in V$ --- see \S2.6), there are only finitely many quaternion division algebras that are unramified at all $v \in V$. So, it suffices to show that for a {\it fixed} quaternion division algebra $D$ over $K$ that is unramified at all $v \in V$, the set of $K$-isomorphism classes of universal covers $\tilde{G} = \widetilde{\mathrm{SU}_n(h,D)}$, where $h$ is a nondegenerate $n$-dimensional skew-hermitian form over $D$, having good reduction at all $v \in V$, is finite.

So, in the remainder of this section, we take $D$ to be a quaternion division algebra that is unramified at all $v \in V$, and we fix a presentation $D = \displaystyle{\left( \frac{a,b}{K} \right)}$. As in \S2.5, using condition (*) and Lemma \ref{L:AdeleRecall}, we may (and will) assume that $V$ arises from an affine model ${\rm Spec}~A$ of $K$, where $A$ is a noetherian unique factorization domain, such that ${\rm char}~K^{(v)} \neq 2$ and $v(a) = v(b) = 0$ for all $v \in V$, and that we have the equality
\begin{equation}\label{E:DivUFD1}
\mathbb{I}(K,V) = \mathbb{I}^{\infty}(K,V) \cdot K^{\times}.
\end{equation}
Thus, we are in the situation of \S2.5; in particular, $a$, $b$, and $2$ are units in $A$ (Remark \ref{R:23}).

Next, for each $v \in V$, let
$$
\Lambda_v := \mathcal{O}_{K_v}\cdot 1 \oplus \mathcal{O}_{K_v}\cdot i \oplus \mathcal{O}_{K_v}\cdot j \oplus \mathcal{O}_{K_v}\cdot k \subset D_v = D\otimes_K K_v.
$$
Clearly, $\Lambda_v$ is an $\mathcal{O}_{K_v}$-order in $D_v$ that is stable under the canonical involution, and in fact, since $a, b, 2 \in \mathcal{O}_{K_v}^{\times}$, it is easy to see that $\Lambda_v$ is a \emph{maximal} order in $D_v$ (cf. \cite[Lemma 15.5.3]{Vo}). Consequently, if $D_v$ is the split algebra, then, using the fact that every maximal order of $M_2(K_v)$ is conjugate to $M_2(\mathcal{O}_{K_v})$ (see \cite[Corollary 10.5.5]{Vo}), we can find an isomorphism $\psi \colon D_v \stackrel{\sim}{\longrightarrow} M_2(K_v)$ such that $\psi(\Lambda_v) = M_2(\mathcal{O}_{K_v}).$ On the other hand, if $D_v$ is a division algebra and $\mathcal{A}_v$ is its valuation ring, then $\Lambda_v \subset \mathcal{A}_v$ (this inclusion follows immediately from the definitions), and since $\mathcal{A}_v$ is an $\mathcal{O}_{K_v}$-order, it follows from maximality that $\mathcal{A}_v = \Lambda_v$. Thus, in both the split and division algebra cases, we can take
\begin{equation}\label{E:good-reduction-forms-convention}
\mathcal{A}_v = \Lambda_v
\end{equation}
for our analysis of skew-hermitian forms with good reduction (see Definition \ref{D:GRskew}).

Finally, let $Q_1$ be the {affine} conic $ax^2 + by^2 = 1$ over $K$ and set $L = K(Q_1).$ Since obviously $K(Q_1) = K(Q)$, where $Q$ is the Severi-Brauer conic associated to $D$, we have the splitting $\varphi \colon D \otimes_K L \to M_2(L)$ given by Lemma \ref{L:SplittingSB}. In the arguments given below, we will work with the model of $L$
$$
\mathcal{Q}_1 = {\rm Spec}~B \to {\rm Spec}~A, \ \ \ \text{where} \ \ B = A[x,y]/(ax^2 + by^2 -1)
$$
and the corresponding divisorial set
$$
\tilde{V}(\mathcal{Q}_1) = \tilde{V}_0(\mathcal{Q}_1) \cup \tilde{V}_1(\mathcal{Q}_1)
$$
that were discussed in \S2.5.

\vskip2mm

Having made these preparatory remarks, we now turn to our argument that will reduce Theorem \ref{T-SB} to Theorem \ref{T-SBUnram}. Suppose that $\{\tilde G_i\}_{i\in I}$ is a family of pairwise non-isomorphic universal covers $\tilde G_i = \widetilde{\mathrm{SU}_n(h_i, D)}$, where each $h_i$ is a nondegenerate $n$-dimensional skew-hermitian form over $D$, such that $\tilde G_i$ has good reduction at all $v \in V$. Since similar forms yield $K$-isomorphic groups, it follows that for $i \neq j$, the forms $h_i$ and $h_j$ are not similar, i.e.\ $h_i$ is not equivalent to $\lambda h_j$ for any $\lambda \in K^\times$. We need to show that $I$ is \emph{finite}.







Now, for a form $h_i$ as above and $v \in V$, let $h_{i,v} := h_i \otimes_K K_v$ denote the base change of $h_i$ to $K_v$.  According to Lemma \ref{L-SUGR}, for each $i \in I$ and any $v \in V$, there exists $\lambda_v^{(i)} \in K_v^{\times}$ such that the form $\lambda_v^{(i)} h_{i,v}$ over $D_v = D \otimes_K K_v$ has good reduction at $v$. Moreover, it easily follows from condition (*) that $h_i$ has good reduction at almost all $v \in V$: indeed, if $C_i \in {\rm GL}_n(D)$ is a Gram matrix of $h_i$, then, by condition (*), the finitely many coordinates, with respect to $1, i, j, k$, of the entries of $C_i$ and $C_i^{-1}$ are $v$-integral for almost all $v \in V$. Thus, for all such $v$, we have $C_i \in {\rm GL}_n(\Lambda_v) = {\rm GL}_n(\mathcal{A}_v)$. Therefore,
we can assume that $\lambda_v^{(i)} = 1$ for all but finitely many $v \in V$. Thus, $\lambda^{(i)} := (\lambda_v^{(i)})_{v \in V}$ defines an element of the idele group $\mathbb{I}(K,V).$ Then, in view of
(\ref{E:DivUFD1}), it follows that we can find, for each $i \in I$, an element $\mu_i \in K^{\times}$ such that $\mu_i h_i$ has good reduction at all $v \in V$.

Consequently,
replacing $h_i$ by $\mu_i h_i$, we can assume that we are given a family $\{h_i \}_{i \in I}$ of $n$-dimensional nondegenerate skew-hermitian forms over $D$ such that

\vskip1mm

\noindent $\bullet$ \parbox[t]{16.5cm}{for each $i \in I$, the form $h_i$ has good reduction at all $v \in V$; and}

\vskip1mm

\noindent $\bullet$ \parbox[t]{16.5cm}{for $i, j \in I$, $i \neq j$, the forms $h_i$ and $h_j$ are not similar.}

\vskip2mm

\noindent We need to show that $I$ is finite. To do this, we first note the following.

\begin{lemma}\label{L:FormsGR}
Suppose $h$ is a nondegenerate $n$-dimensional skew-hermitian form over $D$ that has good reduction at all $v \in V$. Let $q = q_{h_{L}}$ be the $2n$-dimensional quadratic form over $L = K(Q_1)$ that corresponds to $h$ under Morita equivalence via the splitting $\varphi$ of Lemma \ref{L:SplittingSB}. Then $q$ has good reduction at all $\tilde{v} \in \tilde{V}(\mathcal{Q}_1).$
\end{lemma}
\begin{proof}
We consider the claim separately for the sets $\tilde{V}_0(\mathcal{Q}_1)$ and  $\tilde{V}_1(\mathcal{Q}_1)$ of horizontal (or geometric) and vertical valuations, respectively.

\vskip1mm

\noindent {\bf The horizontal case.}  Let $\tilde{v} \in \tilde V_0(\mathcal{Q}_1)$. By Lemma \ref{lem24}(iv), the valuation ring $\mathcal{O}_{L, \tilde{v}}$ contains $B$ and $K$ (in particular $2 \in \mathcal{O}_{L, \tilde{v}}^{\times}$), and hence, by Lemma \ref{L:SplittingSB}, we have $\varphi(D) \subset M_2(\mathcal{O}_{L, \tilde{v}}).$ Let $C \in {\rm GL}_n(D)$ be a Gram matrix of $h$. Then, by Proposition~\ref{prop32I}(ii), $q$ has Gram matrix $S_\varphi(C) = \mathbf{J}^{-1}\Phi_\varphi(C)$. Since $\Phi_\varphi(C)$ and $\Phi_\varphi(C)^{-1} = \Phi_\varphi(C^{-1})$ clearly have entries in $\mathcal{O}_{L,\tilde{v}}$, we see that
$S_\varphi(C) \in {\rm GL}_{2n}(\mathcal{O}_{L,\tilde{v}})$. Since $2 \in \mathcal{O}_{L,\tilde{v}}^\times$, it follows that $S_\varphi(C)$
defines a regular quadratic form over $\mathcal{O}_{L,\tilde{v}}$, and therefore $q$ has good reduction at $\tilde{v}$.


\vskip1mm

\noindent {\bf The vertical case.} Let $\tilde{v} \in \tilde{V}_1(\mathcal{Q}_1)$ be the unique place above $v \in V.$ Denote by $L_{\tilde{v}}$ the completion of $L$ with respect to $\tilde{v}$, with valuation ring $\mathcal{O}_{L_{\tilde{v}}}$, so that we have inclusions $K_v \subset L_{\tilde{v}}$ and $\mathcal{O}_{K_v} \subset \mathcal{O}_{L_{\tilde{v}}}$, and also $2 \in \mathcal{O}_{L_{\tilde{v}}}^{\times}.$ Let
$$
\varphi_{\tilde{v}} \colon D \otimes_K L_{\tilde{v}} \to M_2(L_{\tilde{v}})
$$
be the base change of $\varphi \colon D \otimes_K L \to M_2(L)$ to $L_{\tilde{v}}.$ Since $x, y \in B$, it follows from Lemma \ref{lem24}(iv) that $x,y \in \mathcal{O}_{L, \tilde{v}} \subset \mathcal{O}_{L_{\tilde{v}}}$; combining this with our assumption that $a, b \in \mathcal{O}_{K_v}^{\times}$ then yields
\begin{equation}\label{E:vertical-places-good-reduction-argument}
\varphi_{\tilde{v}}(\Lambda_v) \subset M_2(\mathcal{O}_{L_{\tilde{v}}}).
\end{equation}

Since $h$ has good reduction at $v$, it follows that the form $h_v$ over $D_v$ has a Gram matrix
\begin{equation}\label{E:Gram-matrix-form-good-reduction}
C \in {\rm GL}_n(\mathcal{A}_v) = {\rm GL}_n(\Lambda_v)
\end{equation}
(see Definition \ref{D:GRskew} and the comment following it). Then $C$ is also a Gram matrix of $h_{L_{\tilde{v}}} = h_v\otimes_{K_v} L_{\tilde{v}}$.

Now, we observe that $q\otimes_L L_{\tilde{v}}$ is the quadratic form $q_{h_{L_{\tilde{v}}},\varphi_{\tilde{v}}}$ constructed in
Proposition~\ref{prop32I}(ii), since clearly $\Phi_{\varphi_{\tilde{v}}}(C_0) = \Phi_\varphi(C_0)$ for any Gram matrix $C_0 \in {\rm GL}_n(D)$ of $h$. According to Proposition~\ref{prop32I}(iii), this form does not depend, up to equivalence, on the choice of the Gram matrix of $h_{L_{\tilde{v}}}$. Hence, $q\otimes_L L_{\tilde{v}}$ is equivalent to the quadratic form with Gram matrix $S_{\varphi_{\tilde{v}}}(C) = \mathbf{J}^{-1}\Phi_{\varphi_{\tilde{v}}}(C)$. It follows from (\ref{E:vertical-places-good-reduction-argument}) and (\ref{E:Gram-matrix-form-good-reduction}) that $S_{\varphi_{\tilde{v}}}(C) \in {\rm GL}_{2n}(\mathcal{O}_{L_{\tilde{v}}})$.
Thus, since $2 \in \mathcal{O}_{L_{\tilde{v}}}^\times$, we see that $S_{\varphi_{\tilde{v}}}(C)$
defines a regular quadratic form $\tilde{q}$ over $\mathcal{O}_{L_{\tilde{v}}}$ with $\tilde{q} \otimes_{\mathcal{O}_{L_{\tilde{v}}}} L_{\tilde{v}} \simeq q\otimes_L L_{\tilde{v}}$. We conclude that $q$ has good reduction at $\tilde{v}$.
\end{proof}

Returning to the argument, we see that the quadratic forms $\{q_{h_i}\}_{i \in I}$ over $L= K(Q_1)$ all have good reduction at all $\tilde{v} \in \tilde{V}(\mathcal{Q}_1).$ Suppose we know that there are finitely many equivalence classes of $2n$-dimensional nondegenerate quadratic forms over $L$ that have good reduction at all $\tilde{v} \in \tilde{V}(\mathcal{Q}_1).$ Then the forms $\{q_{h_i}\}_{i \in I}$ fall into finitely many equivalence classes. Consequently, by Proposition~\ref{prop32I}(iv) (applied with $F' = L$ and $\rho = \varphi$),
the skew-hermitian forms $\{h_{i, L}\}_{i \in I}$ fall into finitely many equivalence classes. But then it follows from Proposition \ref{P-PSS} and Witt's cancellation theorem for skew-hermitian forms over $(D,\iota)$ (see, e.g., \cite[Ch.~7, \S9]{Scharlau})
that the forms $\{h_i \}_{i \in I}$ fall into finitely many equivalence classes.
From this, we conclude that the set $I$ is finite, as needed.

Thus, to complete the argument, we need to show that for any $n \geq 2$, there are finitely many equivalence classes of $n$-dimensional nondegenerate quadratic forms over $L$ that have good reduction at all $\tilde{v} \in \tilde{V}(\mathcal{Q}_1).$ According to Proposition \ref{P:FormsGoodReduction}, for this, it suffices to prove that the unramified cohomology groups $H^i(L, \mu_2)_{\tilde{V}(\mathcal{Q}_1)}$ are finite for all $i \geq 1.$ This reduces Theorem \ref{T-SB} to Theorem \ref{T-SBUnram}, which will be proved in the next section.

\vskip3mm

To conclude this section, assuming Theorem \ref{T-SB}, let us now give

\vskip2mm

\noindent {\it Proof of Theorem \ref{T-GenusSpinor}.} Suppose $q$ is a nondegenerate quadratic form over $K = k(C)$ in $n \geq 10$ variables and let $G = {\rm Spin}_n(q)$ be the corresponding spinor group. We note that if $n$ is odd or if $q$ is an isotropic quadratic form of even dimension $\geq 10$, then the finiteness of $\gen_K(G)$ was established in \cite[Theorem 1.2]{Spinor}. Quite generally,  it follows from \cite[Proposition 5.3]{Spinor} that for a nondegenerate quadratic form $q$ over $K$ of even dimension $\geq 10$, any group in $\gen_K(G)$ is either of the form ${\rm Spin}_n(q')$ for an $n$-dimensional nondegenerate quadratic form $q'$ or is the universal cover of a quaternionic special unitary group as considered in Theorem \ref{T-SB}.

Now let $V$ be a divisorial set of places of $K$. After possibly removing finitely many places from $V$ and then picking a divisorial set of places in the complement (see \S2.1), we can assume that ${\rm char}~K^{(v)} \neq 2$ for all $v \in V$ and also that $G$ has good reduction at all $v \in V$. Then by \cite[Theorem 1.1]{CRR-MaxTori}, every $G' \in \gen_K(G)$ also has good reduction at all $v \in V$. According to \cite[Theorem 1.1]{Spinor}, the set of $K$-isomorphism classes of spinor groups of nondegenerate quadratic forms in $n$ variables over $K$ that have good reduction at all $v \in V$ is finite; Theorem \ref{T-SB} provides the analogous statement for the universal covers of special unitary groups. The finiteness of $\gen_K(G)$ now follows. \hfill $\Box$

\section{Proof of Theorem \ref{T-SBUnram}}\label{S:FinUnram}

In this section, we turn to the proof of Theorem \ref{T-SBUnram}, which we restate here for the reader's convenience.

\begin{thm}\label{T-SBUnramRestatement}
Let $K = k(C)$ be the function field of a smooth geometrically integral curve $C$ over a number field $k$.
Suppose $$Q: ax^2 + by^2 = z^2$$ (with $a, b \in K^{\times}$) is an anisotropic conic in $\mathbb{P}^2_K$,
and
consider the function field $L = K(Q).$ Then for any divisorial set of places $V'$ of $L$, the
unramified cohomology groups $H^i(L, \mu_2)_{V'}$ are finite for all $i \geq 1.$
\end{thm}

For the proof, we will consider separately the cases (i) $i = 1, 2$, and $4$; (ii) $i \geq 5$; and (iii) $i = 3.$

\vskip2mm

\noindent {\bf Finiteness of $H^i(L, \mu_2)_{V'}$ for $i = 1, 2$, and $4$.}  For any divisorial set $V'$ of $L$, the finiteness of $H^1(L, \mu_2)_{V'}$ is well-known (see the discussion following \cite[Proposition 4.1]{Spinor}), and the finiteness of $H^2(L, \mu_2)_{V'} = {}_2\Br(L)_{V'}$
is a consequence of the main results of \cite{CRR3}. Furthermore, the fact that $H^4(L, \mu_2)_{V'}$ is finite follows from \cite[Th\'eor\`eme 6.2]{CTFinitude}.

\vskip4mm

\noindent {\bf Finiteness of $H^i(L, \mu_2)_{V'}$ for $i \geq 5.$} The required finiteness in this case is a consequence of the following general statement. To fix notations, given a smooth geometrically integral variety $X$ over a field $\ell$, we denote by $V_0(X)$ the set of geometric places of the function field $F = \ell(X)$, i.e. the set of discrete valuations of $F$ associated with the prime divisors of $X$.

\begin{prop}\label{P-HigherUnramFinite}
Let $F = k(X)$ be the function field of a smooth geometrically integral variety $X$ of dimension $d$ over a number field $k$.
Then the unramified cohomology groups $H^i(F, \mu_2)_{V_0(X)}$ are finite for all $i \geq d+3$.
\end{prop}
\begin{proof}
Fix $i \geq d+3$ and let $r = i - (d+2).$ We set $\tilde{F} = F(t_1, \dots, t_r)$, which we will view as the function field of $\tilde{X} = X \times_k \mathbb{A}_k^r.$  Since $\dim \tilde{X} = d + r = i-2$, the cohomological Hasse principle proved by Jannsen (see \cite[Theorem 0.4]{Jan}) implies that the global-to-local map
$$
H^i(\tilde{F}, \mu_2) \to \bigoplus_{v \in V^k} H^i(\tilde{F}_v, \mu_2)
$$
is injective (here, the sum is taken over the set $V^k$ of all places of $k$ and $\tilde{F}_v = k_v(\tilde{X} \times_k k_v)$). This, in turn, yields an embedding
$$
H^i(\tilde{F}, \mu_2)_{V_0(\tilde{X})} \to \bigoplus_{v \in V^k} H^i(\tilde{F}_v, \mu_2)_{V_0(\tilde{X}_v)}
$$
of unramified cohomology groups, where we set $\tilde{X}_v = \tilde{X} \times_k k_v.$ On the other hand, according to \cite[Theorem 4.1.5]{CT-SB}, the natural map $H^i(F, \mu_2) \to H^i(\tilde{F}, \mu_2)$ induces an isomorphism
$$
H^i({F}, \mu_2)_{V_0(X)} \simeq H^i(\tilde{F}, \mu_2)_{V_0(\tilde{X})},
$$
and similarly for $F_v$. Thus, we obtain an embedding
$$
H^i(F, \mu_2)_{V_0(X)} \hookrightarrow \bigoplus_{v \in V^k} H^i(F_v, \mu_2)_{V_0(X_v)},
$$
where $X_v = X \times_k k_v$. Now, if $v$ is either a finite place or a complex place of $k$, then $H^i (F_v, \mu_2) = 0$ for all $i \geq d+3$ by cohomological dimension considerations. On the other hand, if $v$ is a real place, then $H^i(F_v, \mu_2)_{V_0(X_v)}$ is finite by \cite[Main Theorem]{CTP-Real}, and the finiteness of $H^i(F, \mu_2)_{V_0(X)}$ follows.
\end{proof}

\vskip2mm

\noindent Now, we can view $L$ as the function field $k(X)$ of a smooth (affine) geometrically integral variety $X$ over the number field $k$. By spreading out, we can find a flat normal model $\mathcal{X} \to {\rm Spec}(\mathcal{O}_{k,S})$ of $X$ over a ring $\mathcal{O}_{k,S}$ of $S$-integers (such that the fibers are smooth geometrically integral varieties over the residue fields). The corresponding divisorial set $\tilde{V}(\mathcal{X})$ contains $V_0(X)$, and hence, applying Proposition \ref{P-HigherUnramFinite} with $d=2$, we obtain the finiteness of $H^i(L, \mu_2)_{\tilde{V}(\mathcal{X})}$ for $i \geq 5$. Since any two divisorial sets are commensurable (see \S2.1), it is easy to see that a divisorial set $V'$ of $L$ contains $\tilde{V}(\tilde{\mathcal{X}})$ for some dense open subset $\tilde{X} \subset X$ and its model $\tilde{\mathcal{X}}$. As before, $\tilde{V}(\tilde{\mathcal{X}})$ contains the set of geometric places $V_0(\tilde{X})$, so, applying Proposition \ref{P-HigherUnramFinite} to $\tilde{X}$ yields the finiteness of $H^i(L, \mu_2)_{\tilde{V}(\tilde{\mathcal{X}})}$ for $i \geq 5$ and hence also the finiteness of $H^i(L, \mu_2)_{V'}$.


\vskip4mm

\noindent {\bf Finiteness of $H^3(L, \mu_2)_{V'}$.} We first prove the finiteness of $H^3(L, \mu_2)_{\tilde{V}(\mathcal{Q})}$ for the divisorial set
$\tilde{V}(\mathcal{Q})$ associated with the model $\mathcal{Q}$ of $L$ discussed in \S2.5. The general case then easily follows via an application of
Proposition \ref{P:Open}.


Our setup is as follows. Suppose $D = \displaystyle{\left( \frac{a,b}{K} \right)}$ is a quaternion division algebra over $K$ with corresponding (anisotropic) Severi-Brauer conic $Q \subset \mathbb{P}_K^2$ given by the equation $ax^2 + by^2 = z^2$. Let $V$ be a divisorial set of places of $K$. As in \S2.5, we may assume that $V$ arises from an affine model ${\rm Spec}~A$ of $K$, where $A$ is a noetherian unique factorization domain, such that ${\rm char}~K^{(v)} \neq 2$ and $v(a) = v(b) = 0$ for all $v \in V$. We note that these assumptions imply, in particular, that $D$ is unramified at all $v \in V$.


\vskip1mm

Next, we consider the model of $L = K(Q)$
$$
\mathcal{Q} \to {\rm Spec}~A, \ \ \ \text{where} \ \ \ \mathcal{Q} = \operatorname{Proj}\big(A[x,y,z]/(ax^2+by^2-z^2)\big)
$$
and the corresponding divisorial set of places
$$
\tilde{V}(\mathcal{Q}) = \tilde{V}_0(\mathcal{Q}) \cup \tilde{V}_1 (\mathcal{Q})
$$
that were discussed in \S2.5 (we recall that the elements of $\tilde{V}_0(\mathcal{Q})$ correspond to the closed points of $Q$, whereas the elements of $\tilde{V}_1(\mathcal{Q})$ are associated with the curves $Q^{(v)}$ for $v \in V$, and for each $v \in V$, there is a \emph{unique} $\tilde{v} \in \tilde{V}_1(\mathcal{Q})$ above $v$).
With these notations, we have the following.

\begin{prop}\label{P:UnramDeg3}
The unramified cohomology group $H^3(L, \mu_2)_{\tilde{V}({\mathcal{Q}})}$ is finite.
\end{prop}

We begin with several preliminary observations. First, a well-known consequence of the truth of
the Bloch-Kato Conjecture (proved by Voevodsky) is that for any field $F$ and any prime $p \neq {\rm char}~F$, the map
$$
H^{j}(F, \mu_{p^n}^{\otimes (j -1)}) \stackrel{\iota_{j, n}}{\longrightarrow} H^{j}(F, \Q_p/ \Z_p (j -1))
$$
induced by the natural inclusion $\mu_{p^n}^{\otimes (j -1)} \hookrightarrow \Q_p/ \Z_p (j -1)$
is injective for all $n \geq 1$ and all $j \geq 2$ (cf. \cite[p. 3]{Jan}). In particular, in our case, it is enough to show that
$$
\iota_{3,1}(H^3(L, \mu_2)_{\tilde{V}({\mathcal{Q}})}) \subset H^3(L, \Q_2/ \Z_2(2))
$$
is finite.

Next, suppose $D$ is a quaternion division algebra over a field $F$ of characteristic $\neq 2$ with corresponding Severi-Brauer conic $Q \subset \mathbb{P}^2_F$. Denote by $V_0$ the set of geometric places of the function field $F(Q)$ (i.e. the discrete valuations corresponding to the closed points of $Q$). We have the following important result.

\begin{prop}\label{P:Kahn} {\rm (cf. \cite[Proposition A.1]{KRS})}
The natural map
$$
\eta_F \colon H^3(F, \Q_2 /\Z_2 (2)) \to H^3(F(Q), \Q_2/\Z_2(2))
$$
induces a surjection
$$
\eta_F \colon H^3(F, \Q_2 /\Z_2 (2)) \to H^3(F(Q), \Q_2/\Z_2(2))_{V_0}
$$

\end{prop}

We will also need the next elementary statement.

\begin{lemma}\label{L-KerBaseChange}
For any prime $p \neq {\rm char}~F$, the kernel of the natural map
$$
\eta_F \colon H^3(F, \Q_p/ \Z_p(2)) \to H^3(F(Q), \Q_p/ \Z_p(2))
$$
is annihilated by 2.
\end{lemma}
\begin{proof}
Let $F'/F$ be a (separable) quadratic extension that splits $D$, so that $Q \times_F F' \simeq \mathbb{P}_{F'}^1$.  Then we have the following commutative diagram
$$
\xymatrix{H^3(F, \Q_p/\Z_p(2)) \ar[r]^{\eta_F} \ar[d]_{\gamma} & H^3(F(Q), \Q_p/ \Z_p(2)) \ar[d]^{\beta} \\ H^3(F', \Q_p/ \Z_p (2)) \ar[r]^{\eta_{F'}} & H^3(F'(Q), \Q_p/ \Z_p(2))}
$$
where
$\gamma$ and $\beta$ are restriction maps. Now, $F'(Q) = F'(x)$ is a purely transcendental extension in one variable and it is well known that $\eta_{F'}$ is injective in this case (cf. \cite[Ch. III, \S9]{GMS}).
Thus, if $x \in \ker \eta_F$, we have $(\eta_{F'} \circ \gamma) (x) = (\beta \circ \eta_F)(x)= 0$, hence $\gamma(x) = 0.$ But then, applying corestriction, we conclude that $2x = 0$, as needed.
\end{proof}
Now, identifying $H^3(F, \mu_4^{\otimes 2})$ and $H^3(F(Q), \mu_2)$ with their images under the inclusions
$$
H^3(F, \mu_4^{\otimes 2}) \hookrightarrow H^3(F, \Q_2/ \Z_2(2)) \ \ \ \text{and} \ \ \ H^3(F(Q), \mu_2) \hookrightarrow H^3(F(Q), \Q_2/ \Z_2(2))
$$
discussed above and using Proposition \ref{P:Kahn} in combination with Lemma \ref{L-KerBaseChange}, we obtain
\begin{cor}\label{C-KRS}
The natural map $\eta_F \colon H^3(F, \Q_2/ \Z_2(2)) \to H^3(F(Q), \Q_2/ \Z_2(2))$ induces a surjection
$$
H^3(F, \mu_4^{\otimes 2}) \to H^3(F(Q), \mu_2)_{V_0}.
$$
\end{cor}
\begin{proof}
Given $x \in H^3(F(Q), \mu_2)_{V_0}$, according to Proposition \ref{P:Kahn}, there exists $y \in H^3(F, \Q_2/ \Z_2(2))$ such that $\eta_F(y) = x.$ Since $x$ has order dividing 2, it follows from Lemma \ref{L-KerBaseChange} that $y$ has order dividing 4, hence lies in $H^3(F, \mu_4^{\otimes 2}).$
\end{proof}

\vskip2mm

\noindent {\it Proof of Proposition \ref{P:UnramDeg3}.} We will show that the natural map
$$
\eta_K \colon H^3(K, \Q_2/ \Z_2(2)) \to H^3(L, \Q_2/ \Z_2(2))
$$
induces a surjection
$$
H^3(K, \mu_4^{\otimes 2})_V \to H^3(L, \mu_2)_{\tilde{V}({\mathcal{Q}})}.
$$
Since $H^3(K, \mu_4^{\otimes 2})_V$ is finite by \cite[Theorem 6.1]{Spinor}, this will complete the argument.

So, suppose $t \in H^3(L, \mu_2)_{\tilde{V}({\mathcal{Q}})}$ and let $\xi \in H^3(K, \mu_4^{\otimes 2})$ be such that $\eta_K (\xi) = t$ (which is possible by Corollary \ref{C-KRS}). Fix $\tilde{v} \in \tilde{V}_1(\mathcal{Q})$ and take $v \in V$ such that $\tilde{v} \vert v.$ Since $v(a) = v(b) = 0$  and $\operatorname{char}K^{(v)} \neq 2$, we obtain from Lemma \ref{lem24}(iii) that
the ramification index $e(\tilde{v}|v)$ is $1$, and also that the residue field $L^{(\tilde{v})}$ is the function field $K^{(v)}(Q^{(v)})$ of the reduction $Q^{(v)}\colon \bar ax^2 + \bar by^2 = z^2$, which is the Severi--Brauer conic of the residue algebra $D^{(v)} = \displaystyle{ \left(\frac{\bar a,\bar b}{K^{(v)}} \right)}$ of $D$ at $v$. According to \cite[Ch. II, Proposition 8.2]{GMS}, the fact that $e(\tilde{v}|v) = 1$ implies that we have
the following commutative diagram
$$
\xymatrix{H^3(K, \mu_4^{\otimes 2}) \ar[rr]^{\eta_K} \ar[d]_{\rho^3_v} & & H^3(L, \Q_2/ \Z_2 (2)) \ar[d]^{\rho^3_{\tilde{v}}} \\ H^2(K^{(v)}, \mu_4) \ar[rr]^{\tilde{\eta}_{K^{(v)}}} & & H^2(L^{(\tilde{v})}, \Q_2/ \Z_2 (1))}
$$
where $\rho^3_v$ and $\rho^3_{\tilde{v}}$ are the residue maps associated with $v$ and $\tilde{v}$, respectively, and $\tilde{\eta}_{K^{(v)}}$ is the natural map. Since $D$ is unramified at $v$, it is well known that $D_v = D\otimes_K K_v$ is split if and only if $D^{(v)}$ is split (see \cite[Ch. II, Proposition 7.7]{GMS} or \cite[\S3]{W2}). Now, if $D_v$ is split, then $Q^{(v)} \simeq \mathbb{P}^1_{K^{(v)}}$, so $L^{(\tilde{v})} = K^{(v)}(x)$ and $\tilde{\eta}_{K^{(v)}}$ is injective. Since $\rho_{\tilde{v}}^3(t) = 0$ by our assumption, it follows that $\rho_v^3(\xi) = 0$, i.e.\ $\xi$ is unramified at $v$.

Thus, we may assume that $D_v$ is a quaternion division algebra over $K_v$. Then $D^{(v)}$ is a quaternion division algebra over $K^{(v)}$, and since $L^{(\tilde{v})} = K^{(v)}(Q^{(v)})$, by Amitsur's theorem, $\ker\tilde\eta_{K^{(v)}}$ is the subgroup $\langle[D^{(v)}]\rangle$ of the Brauer group generated by the class of $D^{(v)}$ (cf. \cite[Theorem 5.4.1]{GS}).

Now, the fact that $V$ is a divisorial set implies that the element $\xi \in H^3(K, \mu_4^{\otimes 2})$ is unramified outside a finite subset of places $S \subset V$. Indeed, $\xi$ arises by inflation from an element of $H^3(K'/K,\mu_4^{\otimes 2})$ for some finite Galois extension $K'/K$. Since $V$ satisfies (*), the extension $K'/K$ is unramified at all places $v \in V\setminus S$ for some finite subset $S \subset V$.
It then follows from the definition of the residue maps that $\xi$ is unramified at all $v \in V \setminus S$ (cf. \cite[Ch. II, \S7]{GMS}).

From this discussion and the fact that $t= \eta_K (\xi)$ is unramified at all $\tilde{v} \in \tilde{V}_1(\mathcal{Q})$, we conclude that
$$
\rho^3_v(\xi) \in  \langle [D^{(v)}] \rangle \ \ \ \text{for all} \ v \in S,
$$
and after removing from $S$ the places where $\rho_v^3(\xi) = 0$, we may assume that $\rho_v^3(\xi) = [D^{(v)}]$ for all $v \in S$. If $S$ becomes empty, then $\xi$ is unramified at all $v \in V$, i.e.\ $\xi \in H^3(K,\mu_4^{\otimes 2})_V$ with $\eta_K(\xi) = t$, and we are done. So, we may assume that $S \neq \emptyset$.
Since, by our construction, $V$ arises from an affine model ${\rm Spec}~A$ where $A$ is a unique factorization domain, we can find $\alpha \in A$ such that
$$
v(\alpha) = \left\{ \begin{array}{ll} 1, & v \in S \\ 0, & v \in V \setminus S. \end{array} \right.
$$
Let $\chi_{\alpha} \in H^1(K, \mu_4) \simeq K^{\times}/{K^{\times}}^4$ be the corresponding cohomology class and set
$$
\zeta = [D] \cup \chi_{\alpha} \in H^3(K, \mu_4^{\otimes 2}).
$$
Then it follows from \cite[Lemma 7.1]{Srimathy} that $\rho^3_v(\xi) = \rho^3_v(\zeta)$ for all $v \in V$. Consequently, $\xi - \zeta \in H^3(K, \mu_4^{\otimes 2})_{V}$. On the other hand, since $D$ splits over $L$, we have
$$
\eta_K(\xi - \zeta) = \eta_K(\xi) = t.
$$
In conclusion, we obtain a surjection
$$
H^3(K, \mu_4^{\otimes 2})_V \to H^3(L, \mu_2)_{\tilde{V}({\mathcal{Q}})},
$$
as needed, completing the proof. \hfill $\Box$

\vskip3mm

Next, let $\tilde Q \subset Q$ be an open subvariety, obtained by removing finitely many closed points $P_j$ \textup{(}$j \in J$\textup{)} from $Q$, and let $Y_j$ be the closure of $P_j$ in $\cQ$. By Lemma~\ref{lem23}(ii), after replacing $\Spec A$ by a principal open subset $\Spec A_f$ (which does not affect the preceding discussion), we may assume that the $Y_j$ are pairwise disjoint, normal, and finite over $\Spec A$.

We now set $\tilde\cQ := \cQ\setminus\bigcup_{j\in J}Y_j$. By Lemma~\ref{lem23}(ii), $\tilde\cQ$ is a model of $L$ over $\Spec A$ with smooth geometrically integral fibers. Let $\tilde V(\tilde\cQ)$ be the corresponding divisorial set of places of $L$ (which is obtained from $\tilde V(\cQ)$ simply by removing the finitely many horizontal places corresponding to the points $P_j$).

By construction, $Y := \cQ\setminus\tilde\cQ = \bigcup_{j\in J}Y_j$ is a normal closed subscheme of $\cQ$ of pure codimension $1$, whose irreducible components are the horizontal prime divisors $Y_j$. Let $\kappa_j$ be the function field of $Y_j$, i.e.\ the residue field of $P_j$. Then $\kappa_j$ is a finite extension of $K$, hence a $2$-dimensional global field, and $Y_j$ is a model of $\kappa_j$. Consequently, according to the main results of \cite{CRR3}, the unramified cohomology groups $H^2(\kappa_j,\mu_2)_{V(Y_j)}$ are finite. Since $\cQ$ is an excellent (being of finite type over $\mathbb{Z}$) integral noetherian normal scheme with function field $L$ on which $2$ is invertible, Proposition \ref{P:Open} (applied to the case
$X = \cQ$, $U = \tilde\cQ$, $i = 3$, and $n = 2$) yields an exact sequence
\[
0 \to H^3(L,\mu_2)_{\tilde V(\cQ)} \to H^3(L,\mu_2)_{\tilde V(\tilde\cQ)} \to \bigoplus_{j\in J} H^2(\kappa_j,\mu_2)_{V(Y_j)}.
\]
Since the left term is finite by Proposition \ref{P:UnramDeg3} and, as just noted, the right term is finite, it follows that
the group $H^3(L,\mu_2)_{\tilde V(\tilde\cQ)}$ is finite.

Finally, let $V'$ be an arbitrary divisorial set of places of $L$. Since any two divisorial sets are commensurable (see \S2.1), $V' \supset \tilde V(\cQ)\setminus S'$ for a finite subset $S' \subset \tilde V(\cQ)$. The vertical places in $S'$ lie over finitely many $v \in V$, and we can remove them by replacing $\Spec A$ by a principal open subset. The horizontal places in $S'$ correspond to finitely many closed points $P_j$ of $Q$, which we remove from $Q$ to obtain $\tilde Q$. Then $V' \supset \tilde V(\tilde\cQ)$, so $H^3(L,\mu_2)_{V'} \subset H^3(L,\mu_2)_{\tilde V(\tilde\cQ)}$ is finite. This completes our proof of the finiteness of the degree $3$ unramified cohomology of $L$.

\section{Proof of Theorem \ref{T-Surface}}\label{S:Surface}

In this section, we will give the proof of Theorem \ref{T-Surface}, which, we recall, deals with
the finiteness properties of certain groups with good reduction over function fields of $p$-adic surfaces. This result will be deduced from Proposition \ref{P-FinUnramPadic} below, which establishes the finiteness of some unramified cohomology groups of such fields.


Our notations throughout this section will be as follows. Let $k$ be a $p$-adic field and $S$ a smooth geometrically integral surface over $k$. We set $K = k(S)$ to be the function field of $S$ and take $V$ to be the set of geometric places of $K$ corresponding to the prime divisors of $S$. 

\begin{prop}\label{P-FinUnramPadic}
In the above notations, we
have the following:

\vskip1mm

\begin{enumerate}[label=\textup{(\roman*)}, leftmargin=*]
  \item The unramified cohomology groups $H^i(K,\mu_2)_V$ are finite for $i = 1, 2,$ and $4$.
  \item If $p$ is odd, then $H^3(K, \mu_2)_V$ is finite.
\end{enumerate}

\end{prop}

\vskip1mm

\noindent (We note that by cohomological dimension considerations, we have $H^i(K, \mu_2) = 0$ for $i \geq 5$ --- cf. \cite[Ch. II, \S\S2 and 4]{Serre-GC}.)

\vskip2mm

For the argument, we will need the following well-known statement, whose proof uses the finiteness of \'etale cohomology groups over an algebraically closed field, finiteness properties of fields of type (F),
and the Hochschild-Serre spectral sequence --- see, e.g., \cite[Corollary 3.2]{IRap} for the details and a generalization.

\begin{thm}\label{T-EtaleFinite}
Suppose $X$ is a smooth geometrically integral variety over a field $k$ that is of characteristic 0 and of type ${\rm (F)}.$
Then for any positive integer $m$, the \'etale cohomology groups $\he^i(X, \mu_{m,X}^{\otimes j})$ are finite for all $i \geq 0$ and all $j$.
\end{thm}

\vskip2mm

\noindent {\it Proof of Proposition \ref{P-FinUnramPadic}.} (i) By cohomological purity for smooth varieties over fields, we have
surjections
$$
\he^1(S, \mu_{2, S}) \twoheadrightarrow H^1(K, \mu_2)_{V} \ \ \ \text{and} \ \ \ \he^2(S, \mu_{2, S}) \twoheadrightarrow H^2(K, \mu_2)_{V}
$$
(see \cite[\S3.4]{CT-SB}). It is well-known that $k$ is of type (F). Thus, Theorem \ref{T-EtaleFinite} yields the finiteness of $H^1(K, \mu_2)_{V}$ and $H^2(K, \mu_2)_{V}$.
The finiteness of $H^4(K, \mu_2)_{V}$ was established in \cite[Th\'eor\`eme 5.1]{CTFinitude}.


\vskip2mm

\noindent (ii) We now assume that $p$ is an {\it odd} prime.
When $S$ is {\it projective}, the finiteness of $H^3(k(S), \mu_2)_V$ was established in \cite[Corollary 2.6]{SaSu} using the Merkurjev-Suslin theorem and considerations from higher-dimensional class field theory.
Let us now consider the case of an arbitrary smooth geometrically integral surface $S$ over $k$ (where the argument is similar to the proof of \cite[Th\'eor\`eme 5.1]{CTFinitude}). First, replacing $S$ by an affine open subset (which only shrinks the set $V$), we may suppose that $S$ is affine. Then taking the Zariski closure of $S$ in a suitable projective space and applying Hironaka's theorem on resolution of singularities, we can realize $S$ as an open subset of a smooth geometrically integral projective surface $\tilde{S}$ over $k$. Let $Y = \tilde{S} \setminus S$ be the complement of $S$, endowed with the structure of a reduced subscheme of $\tilde{S}$. Possibly enlarging $S$, we may assume that $Y$ has pure dimension 1. Let $F = Y_{\rm sing}$ be the singular locus of $Y$, which consists of a finite number of closed points of $\tilde{S}$. We can write
$$
Y \setminus F = \bigcup_{i \in I} Z_i,
$$
where the $Z_i$ are smooth curves over $k$, and we denote by $\kappa_i$ the function field of $Z_i$. Furthermore, we set $X = \tilde{S} \setminus F$. Since $F$ is a finite set of closed points, it has codimension 2 in $\tilde{S}$, so $X$ and $\tilde{S}$ have the same prime divisors. Thus, $V_X = V_{\tilde{S}}$, where $V_X$ and $V_{\tilde{S}}$ are the sets of discrete valuations of $K$ corresponding to the prime divisors of $X$ and $\tilde{S}$, respectively; in particular, $H^3(K, \mu_2)_{V_X} = H^3(K, \mu_2)_{V_{\tilde{S}}}$ is finite by the projective case treated above. Then by Proposition \ref{P:Open},
we have an exact sequence of unramified cohomology groups
$$
H^3(K, \mu_2)_{V_X} \to H^3(K, \mu_2)_V \to \bigoplus_{i \in I} H^2(\kappa_i, \mu_2)_{V_{Z_i}}.
$$
Thus, to obtain the finiteness of $H^3(K, \mu_2)_V$, it suffices to note that, for all $i \in I$, we have surjections
$$
\he^2(Z_i, \mu_{2, Z_i}) \twoheadrightarrow H^2(\kappa_i, \mu_2)_{V_{Z_i}}
$$
by cohomological purity, so all of the groups $H^2(\kappa_i, \mu_2)_{V_{Z_i}}$ are finite by Theorem \ref{T-EtaleFinite}. $\Box$

\vskip3mm

\noindent {\it Proof of Theorem \ref{T-Surface}.}

\vskip2mm

\noindent (i) In view of Proposition \ref{P-FinUnramPadic} and \cite[Theorem 2.1]{Spinor}, it suffices to establish the finiteness of the quotient $\Pic(K, V)/2 \cdot \Pic(K,V)$, where $\Pic(K, V)$ is the Picard group of $V$ (cf. \S\ref{S:Recollections}). For this, we first note that it follows from the definitions that  $\mathrm{Pic}(K,V)$ coincides with the usual Picard group $\mathrm{Pic}(S)$ of $S$. Now, we consider the following Kummer sequence of \'etale sheaves on $S$
$$
1 \to \mu_{2, S} \to \mathbb{G}_{m,S} \stackrel{[2]}{\longrightarrow} \mathbb{G}_{m, S} \to 1,
$$
which gives rise to the sequence
$$
\he^1(S, \mathbb{G}_{m,S}) \stackrel{[2]}{\longrightarrow} \he^1(S, \mathbb{G}_{m,S}) \to \he^2(S, \mu_{2, S})
$$
in \'etale cohomology. Since $\Pic(S) = \he^1(S, \mathbb{G}_{m,S})$ (cf. \cite[Ch. III, Proposition 4.9]{Milne}), we thus obtain an injection
$$
\Pic(K,V)/ 2 \cdot \Pic(K,V) \hookrightarrow \he^2(S, \mu_{2,S}),
$$
and the finiteness of $\Pic(K,V)/ 2 \cdot \Pic(K,V)$ follows from Theorem \ref{T-EtaleFinite}.

\vskip2mm

\noindent (ii) By the proof of \cite[Theorem 8.1]{Spinor}, our statement for special unitary groups associated with quadratic extensions
follows from part (i).

\vskip2mm

\noindent (iii) As shown in the proof of \cite[Theorem 9.1]{Spinor}, our claim follows from the finiteness of $H^3(K, \mu_2)_V$, which we established in Proposition \ref{P-FinUnramPadic} above.

\vskip2mm

\noindent (iv) Since $H^5(K, \mu_2) = 0$ by cohomological dimension considerations, the proof of \cite[Theorem 1.12]{CRR-MaxTori} shows that the required finiteness statement follows from the finiteness of $H^3(K, \mu_2)_V.$

\section{Consequences for the global-to-local map}\label{S:LGConsequence}

As we already mentioned in \S\ref{S:Intro}, if $G$ is a (connected) reductive group over a finitely generated field $K$ and $V$ is a divisorial set of places of $K$, then the expectation is that the global-to-local map
$$
\lambda_{G,V} \colon H^1(K, G) \to \prod_{v \in V} H^1(K_v, G)
$$
is {\it proper} (see \cite[Conjecture 6.1]{RR-survey}). One also expects a similar properness result to hold in the geometric setting of function fields of varieties over fields of characteristic 0 and type (F) equipped with their geometric sets of places.

In \cite[\S6]{RR-survey}, we observed that the truth of (AFC) and (GFC) for an absolutely almost simple simply connected group $G$ would yield the Properness Conjecture for the corresponding adjoint group $\overline{G}.$ In particular, (AFC) and (GFC) directly imply the Properness Conjecture for groups of type $\mathsf{F}_4$ and $\mathsf{G}_2$. On the other hand, our results in \cite{Spinor} and \cite{CRR-MaxTori} established a connection between the finiteness properties of certain unramified cohomology groups and (AFC) and (GFC) for groups of type $\mathsf{G}_2$ as well as groups of type $\mathsf{F}_4$ that split over a quadratic extension of the field of definition (see \cite[\S9]{Spinor} and \cite[Theorem A2.6]{CRR-MaxTori}, respectively). Furthermore, in \cite[Theorem 3.4]{Spinor}, we observed that the finiteness of unramified cohomology groups with $\mu_2$-coefficients implies the Properness Conjecture for special orthogonal groups of quadratic forms (a similar result for special unitary groups associated with quadratic extensions was also deduced in \cite[Theorem 8.3]{Spinor}).

With these facts in mind, we now indicate several consequences of Theorem \ref{T-SBUnram} and Proposition \ref{P-FinUnramPadic} for the Properness Conjecture. In the statement below, we will assume that $(F,V)$ is one of the following:

\vskip2mm

\noindent $\bullet$ \parbox[t]{16.5cm}{$F = K(Q)$, where $K = k(C)$ is a 2-dimensional global field of characteristic 0, $Q: ax^2 + by^2 = z^2$ (with $a, b \in K^{\times}$) is an anisotropic conic in $\mathbb{P}^2_K$, and $V$ is a divisorial set of places of $F$; or}

\vskip2mm

\noindent $\bullet$ \parbox[t]{16.5cm}{$F = k(S)$, where $k$ is a $p$-adic field, with $p$ an odd prime, $S$ is a smooth geometrically integral surface over $k$, and $V$ is the set of geometric places of $F$.}

\vskip2mm

\begin{thm}\label{T:LocGlConsequence}
Consider a pair $(F,V)$ as above. Then the global-to-local map
$$
\lambda_{G,V} \colon H^1(F, G) \to \prod_{v \in V} H^1(F_v, G)
$$
is proper in the following cases:
\vskip1mm
\noindent $\bullet$ \parbox[t]{16cm}{$G = {\rm SO}_n(q)$, where $q$ is a nondegenerate quadratic form over $F$ in $n \geq 5$ variables;}

\vskip1mm

\noindent $\bullet$ \parbox[t]{16cm}{$G= {\rm SU}_n(L/F, h)$, where $L/F$ is a quadratic extension and $h$ is a nondegenerate hermitian $L/F$-form in $n \geq 2$ variables;}

\vskip1mm

\noindent $\bullet$ $G$ is a simple algebraic group of type $\mathsf{G}_2$ over $F$;

\vskip1mm

\noindent $\bullet$ \parbox[t]{16cm}{$G$ is a simple algebraic group of type $\mathsf{F}_4$ over $F$ that splits over a quadratic extension of $F$.}

\end{thm}

\vskip2mm

For some other recent developments on local-global principles in the context of semi-global fields (such as functions fields of curves over $p$-adic fields), we refer the reader to \cite{CTHHKPS}, \cite{CTHHKPS1}, and \cite{PS}.

\section{Proofs of Theorems \ref{T-GenusHilbert} and \ref{T-GenusSurface}}\label{S:GenusHilbert}

Our goal in this section is to prove Theorem \ref{T-GenusHilbert}, and then deduce Theorem \ref{T-GenusSurface}. Since the argument for Theorem \ref{T-GenusHilbert} makes extensive use of generic tori, we begin by quickly recalling the main definitions --- we refer the reader to \cite{PR-Irred}, \cite{PR-MRL1}, and \cite{PR-Vin} for general constructions of generic tori and to \cite{CRR-MaxTori} for applications of generic tori to the analysis of good reduction for groups over finitely generated fields.

Given an algebraic torus $T$ defined over a field $k$, we denote by $k_T$ the minimal splitting field of $T$.
It is well-known that the Galois group $\mathscr{G}_T = \Ga(k_T/k)$ acts faithfully on the group of characters $X(T)$. Now suppose that $G$ is a semisimple $k$-group, $T$ a maximal $k$-torus of $G$, and $\Phi(G , T)$
the corresponding root system. Then the action of $\mathscr{G}_T$ on $X(T)$ permutes the roots, yielding a group homomorphism
$$
\theta_T \colon \mathscr{G}_T \longrightarrow \mathrm{Aut}(\Phi(G ,
T)) \ \ (\subset \mathrm{GL}(X(T) \otimes_{\mathbb{Z}} \mathbb{Q})).
$$
Since $\Phi(G , T)$ generates a finite index subgroup of $X(T)$, this homomorphism is {injective}.
We say that $T$ is {\it generic} over $k$, or $k$-{\it generic}, if the image of $\theta_T$ contains the Weyl group $W(G , T)$.

For the proof of Theorem \ref{T-GenusHilbert}, we will need the following existence result for generic tori over Hilbertian fields (for background on such fields, the reader can consult \cite[Ch. 12, 13]{FrJar} or \cite{SerreGT}).

\begin{prop}\label{P-GenericHilbert}
Let $G$ be a (connected) semisimple algebraic group over a Hilbertian field $k$ of characteristic 0. Then for any finite field extension $\ell/k$, the group $G$ contains a maximal $k$-torus that is generic over $\ell.$
\end{prop}

The proof mimics the argument given by V.E.~Voskresenski\={\i} for semisimple algebraic groups over number fields (see \cite[4.2]{Voskr}). More precisely, the part of the argument dealing with the variety of maximal tori and yielding the existence of a generic torus over a purely transcendental extension of the base field does not require any changes at all. On the other hand, the
specialization part that produces a desired generic torus defined over the base field relies solely on the Hilbertian property of number fields,
which holds in our case by assumption. Thus, to avoid repetition, we will omit the details of the proof here.

\vskip4mm

In the remainder of this section, we will establish the next result, from which Theorem \ref{T-GenusHilbert} will easily follow.

\begin{thm}\label{T:GR-all}
Let $\mathcal{K}$ be a field complete with respect to a discrete valuation $v$ such that the residue field $k$ is a Hilbertian field of characteristic 0, and let $G$ be an absolutely almost simple $\mathcal{K}$-group that has good reduction at $v$. 
Then any  $G' \in \gen_{\mathcal{K}}(G)$  also has good reduction at $v$.
\end{thm}

The proof of Theorem \ref{T:GR-all} uses the strategy we developed earlier in \cite[\S6]{CRR-MaxTori}, and, in particular, requires separate consideration of the cases where the type of $G$ is different from $\mathsf{A}_1$ and $\mathsf{B}_{\ell}$ ($\ell \geq 2$), and where it is one of those types. So, we will simply quote the results from \cite{CRR-MaxTori} that can be transferred into our situation either unaltered or with minimal changes, and then give detailed formulations of those statements that rely on Proposition \ref{P-GenericHilbert}, which was not present in \cite{CRR-MaxTori}.


\begin{thm}\label{T:GR-not} {\rm \cite[Theorem 6.2]{CRR-MaxTori}} Let $\mathcal{K}$ be a field complete with respect to a discrete valuation $v$, and let $G$ be an absolutely almost simple algebraic
$\mathcal{K}$-group of type different from $\textsf{A}_1$ and $\textsf{B}_{\ell}$ $(\ell \geq 2)$. Assume that $G$ contains a maximal $\mathcal{K}$-torus $T$ that is $\mathcal{K}$-generic and whose minimal splitting field $\mathcal{K}_T$ is unramified over $\mathcal{K}$. Then $G$ has good reduction at $v$.
\end{thm}

We also have the following partial converse.
\begin{prop}\label{P:tor10}
Let $\mathcal{K}$ be a field complete with respect to a discrete valuation $v$ and assume that the residue field $k = \mathcal{K}^{(v)}$ is a Hilbertian field of characteristic 0. If $G$ is an absolutely almost simple algebraic $\mathcal{K}$-group that has good reduction at $v$, then $G$ possesses a maximal $\mathcal{K}$-torus $T$ that is $\mathcal{K}$-generic and whose minimal splitting field $\mathcal{K}_T$ is unramified over $\mathcal{K}$.
\end{prop}
\begin{proof}
Let $\mathscr{G}$ be a reductive model of $G$ over the valuation ring $\mathcal{O} \subset \mathcal{K}$. Then the reduction $\underline{\mathscr{G}}$ (i.e. the special fiber of $\mathscr{G}$) is an absolutely almost simple algebraic $k$-group of the same type as $G$ (see \cite[Exp. XXII, Proposition 2.8]{SGA3}). Since $k$ is a Hilbertian field of characteristic 0, according to Proposition \ref{P-GenericHilbert}, there exists a maximal $k$-torus $\overline{\mathscr{T}}$ of $\underline{\mathscr{G}}$ that is generic over $k$. Let $\mathscr{T}$ be a lift of $\overline{\mathscr{T}}$ to $\mathscr{G}$ (cf. \cite[Corollary B.3.5]{Conrad}). Then the generic fiber $T$ is a maximal $\mathcal{K}$-torus of $G$ that is $\mathcal{K}$-generic and whose splitting field $\mathcal{K}_T$ is unramified over $\mathcal{K}$.
\end{proof}

We now deduce

\begin{cor}\label{C:notA1Bl}
Let $\mathcal{K}$ be a field complete with respect to a discrete valuation $v$ and $G$ be an absolutely almost simple algebraic $\mathcal{K}$-group of type different from $\textsf{A}_1$ and $\textsf{B}_{\ell}$ $(\ell \geq 2)$. Assume that the residue field $k = \mathcal{K}^{(v)}$ is a Hilbertian field of characteristic 0 and $G$ has good reduction at $v$. Then any $G' \in \gen_{\mathcal{K}}(G)$ has good reduction at $v$.
\end{cor}
\begin{proof}
According to Proposition \ref{P:tor10}, the group $G$ contains a maximal $\mathcal{K}$-torus $T$ that is generic over $\mathcal{K}$ and whose splitting field $\mathcal{K}_T$ is unramified over $\mathcal{K}$. By definition, any $G' \in \gen_{\mathcal{K}}(G)$ contains a maximal $\mathcal{K}$-torus $T'$ isomorphic to $T$. Clearly $\mathcal{K}_{T'} = \mathcal{K}_{T}$ is unramified over $\mathcal{K}$. Furthermore, the fact that $G'$ is an inner twist of $G$ implies that the subgroups $\theta_T(\mathrm{Gal}(\mathcal{K}_T/\mathcal{K})) \subset \mathrm{Aut}(\Phi(G , T))$ and $\theta_{T'}(\mathrm{Gal}(\mathcal{K}_{T'}/\mathcal{K})) \subset \mathrm{Aut}(\Phi(G , T'))$ have isomorphic images in the groups of symmetries of the corresponding Dynkin diagrams. Hence, the fact that $T$ is $\mathcal{K}$-generic implies the same for $T'.$ Consequently,
$G'$ has good reduction by Theorem \ref{T:GR-not}.
\end{proof}

Next, we turn our attention to groups of types $\textsf{A}_1$ and $\textsf{B}_{\ell}$, where we have the following criterion for the existence of good reduction.

\begin{thm}\label{T:AB}
Let $\mathcal{K}$ be a field that is complete with respect to a discrete valuation $v$ with residue field $k$.
Suppose that $G$ is an absolutely almost simple algebraic $\mathcal{K}$-group of type either $\textsf{A}_1$ or $\textsf{B}_{\ell}$ $(\ell \geq 2)$. 

\vskip2mm

\begin{enumerate}[label=\textup{(\roman*)}, leftmargin=*]
  \item Assume that $\mathrm{char}\: k \neq 2$ if $G$ is of type $\textsf{B}_{\ell}$. If $G$ contains two maximal $\mathcal{K}$-tori $T_1$ and $T_2$ that are generic over $\mathcal{K}$ and whose splitting fields
$\mathcal{K}_{T_i}$ are unramified extensions of $\mathcal{K}$ and satisfy $\mathcal{K}_{T_1} \cap \mathcal{K}_{T_2} = \mathcal{K}$, then $G$ has good reduction at $v$.
  \item Conversely, if $k$ is a Hilbertian field of characteristic 0 and $G$ has good reduction at $v$, then $G$ contains two maximal $\mathcal{K}$-tori $T_1$ and $T_2$ that are generic over $\mathcal{K}$ and whose splitting fields $\mathcal{K}_{T_1}$ and $\mathcal{K}_{T_2}$ are unramified extensions of $\mathcal{K}$ and satisfy $\mathcal{K}_{T_1} \cap \mathcal{K}_{T_2} = \mathcal{K}$.
\end{enumerate}
\end{thm}
\begin{proof}
For the argument, we may assume, without loss of generality, that $G$ is simply connected.

\vskip1mm

\noindent (i) This implication is \cite[Theorem 6.6(1)]{CRR-MaxTori}.

\vskip1mm

\noindent (ii) Let $\mathscr{G}$ be a reductive model of $G$ over the valuation ring $\mathcal{O} \subset \mathcal{K}$ with reduction $\underline{\mathscr{G}}$, which is an absolutely almost simple algebraic $k$-group of the same type as $G$. Since $k$ is a Hilbertian field of characteristic 0, by Proposition \ref{P-GenericHilbert},
we can find a maximal $k$-torus $\overline{\mathscr{T}}_1$ of $\underline{\mathscr{G}}$ that is generic over $k$. Next,
we can again apply Proposition \ref{P-GenericHilbert} to find a maximal $k$-torus $\overline{\mathscr{T}}_2$ of $\underline{\mathscr{G}}$ that is generic over the splitting field $k_{\overline{\mathscr{T}}_1}$ of $\overline{\mathscr{T}}_1$. Now, the fact that the Dynkin diagrams of the types $\textsf{A}_1$ and $\textsf{B}_{\ell}$ do not have nontrivial automorphisms implies that the degrees $[k_{\overline{\mathscr{T}}_i} : k]$ for $i = 1, 2$ are equal to the order $w$ of the Weyl group. Besides, the degree $[k_{\overline{\mathscr{T}}_1} k_{\overline{\mathscr{T}}_2} : k_{\overline{\mathscr{T}}_1}]$ also equals $w$. It follows that
\begin{equation}\label{E:k-disj}
k_{\overline{\mathscr{T}}_1} \cap k_{\overline{\mathscr{T}}_2} = k.
\end{equation}
Let $\mathscr{T}_1$ and $\mathscr{T}_2$ be the lifts of $\overline{\mathscr{T}}_1$ and $\overline{\mathscr{T}}_2$ to $\mathscr{G}$, and let $T_1$ and $T_2$ be the corresponding generic fibers. Then $T_1$ and $T_2$ are maximal $\mathcal{K}$-tori of $G$ that are generic over $\mathcal{K}$ and whose splitting fields $\mathcal{K}_{T_i}$ are unramified extensions of $\mathcal{K}$ with residue fields $k_{\overline{\mathscr{T}}_i}$ for $i = 1, 2$. Then (\ref{E:k-disj}) implies that
$\mathcal{K}_{T_1} \cap \mathcal{K}_{T_2} = \mathcal{K}$, as required.

\end{proof}

From this result, we obtain the following analogue of Corollary \ref{C:notA1Bl} for groups of types $\textsf{A}_1$ and $\textsf{B}_{\ell}$.

\begin{cor}\label{C:GR-AB}
Suppose $\mathcal{K}$ is a field that is complete with respect to a discrete valuation $v$ and let $G$ be an absolutely almost simple algebraic $\mathcal{K}$-group of type either $\textsf{A}_1$ or $\textsf{B}_{\ell}$ $(\ell \geq 2)$ that has good reduction at $v$. Assume that the residue field $k = \mathcal{K}^{(v)}$ is a Hilbertian field of characteristic 0. Then any  $G' \in \gen_{\mathcal{K}}(G)$ also has good reduction at $v$.
\end{cor}
\begin{proof}
Since $G$ has good reduction at $v$, according to Theorem \ref{T:AB}(ii), we can find two maximal $\mathcal{K}$-tori $T_1$ and $T_2$ of $G$ that are generic over $\mathcal{K}$ and whose splitting fields are unramified extensions of $\mathcal{K}$ and satisfy $\mathcal{K}_{T_1} \cap \mathcal{K}_{T_2} = \mathcal{K}$. On the other hand, $G'$ contains maximal $\mathcal{K}$-tori $T'_1$ and $T'_2$ that are $\mathcal{K}$-isomorphic to $T_1$ and $T_2$, respectively. Clearly, $T'_1$ and $T'_2$ have properties analogous to those of $T_1$ and $T_2$, so $G'$ has good reduction by Theorem \ref{T:AB}(i).
\end{proof}

We observe that Theorem \ref{T:GR-all} now follows from Corollaries \ref{C:notA1Bl} and \ref{C:GR-AB}. To complete the proof of Theorem \ref{T-GenusHilbert}, we will also need the following well-known

\begin{lemma}\label{L-TrDegHilbert}
Suppose $k$ is a field of characteristic 0 and $F/k$ is a finitely generated field extension of transcendence degree $\geq 1.$ Then $F$ is Hilbertian.
\end{lemma}
\begin{proof}
Let $d = {\rm tr. deg.}_k F.$ Then we can write $F$ as a finite (separable) extension of the rational function field $k' = k(x_1, \dots, x_d).$ The field $k'$ is Hilbertian by \cite[Proposition 13.2.1]{FrJar}, and then $F$ is Hilbertian by \cite[Corollary 12.2.3]{FrJar}.
\end{proof}

\vskip2mm

\noindent {\it Proof of Theorem \ref{T-GenusHilbert}.} Let $v$ be a geometric place of $K = k(X)$, where $X$ is a smooth geometrically integral variety of dimension $> 1$ over a field $k$ of characteristic 0. Let $\mathcal{K} = K_v$ be the corresponding completion and $K^{(v)}$ be the residue field. Note that since ${\rm tr. deg}_k K^{(v)} = \dim X -1 \geq 1$, we see that $K^{(v)}$ is Hilbertian by Lemma \ref{L-TrDegHilbert}. Moreover, according to \cite[Corollary 3.3]{CRR-MaxTori}, for any $G' \in \gen_K (G)$, we have $G' \times_K \mathcal{K} \in \gen_{\mathcal{K}} (G \times_K \mathcal{K}).$ Thus, Theorem \ref{T-GenusHilbert} follows immediately from Theorem \ref{T:GR-all}. $\Box$

\vskip2mm

\noindent {\it Proof of Theorem \ref{T-GenusSurface}.} We note that parts (ii) and (iii) follow immediately from Theorems \ref{T-Surface} and \ref{T-GenusHilbert}. For part (i), one additionally needs to use \cite[Proposition 5.3]{Spinor}, according to which, under our assumptions, any group in the genus of ${\rm Spin}_n(q)$ is of the form ${\rm Spin}_n(q')$. \hfill $\Box$

\section{Proof of Theorem \ref{T-HigherSha}}\label{S:HigherSha}

In this section, we will establish Theorem \ref{T-HigherSha} on the finiteness of the higher Tate-Shafarevich groups of a torus over the function field of a variety defined over a base field that is of characteristic 0 and of type (F).

We begin with the following elementary lemma, which did not appear in the previous results in \cite{RR-Tori}, \cite{RR3}, and \cite{RR-TS} on the finiteness of Tate-Shafarevich groups.


\begin{lemma}\label{HSL:1}
Let $\Gamma$ be a finite group of order $n$ and $\Delta$ be a (commutative) $\Gamma$-module. If the groups ${}_n\Delta = \{ \delta \in \Delta \ \vert \ n\delta = 0 \}$ and $\Delta/n\Delta$ are finite, then the cohomology groups $H^i(\Gamma , \Delta)$ are finite for all $i \geq 1$.
\end{lemma}
\begin{proof}
We have the following two exact sequences of $\Gamma$-modules:
$$
0 \to {}_n\Delta \longrightarrow \Delta \stackrel{[n]}{\longrightarrow} n\Delta \to 0
$$
(where $[n]$ denotes multiplication by $n$), and
$$
0 \to n\Delta \longrightarrow \Delta \longrightarrow \Delta / n\Delta \to 0.
$$
For any $i \geq 1$, the associated sequences of cohomology groups look as follows:
\begin{equation}\label{HSE:1}
H^i(\Gamma , {}_n\Delta) \stackrel{\lambda^i}{\longrightarrow} H^i(\Gamma , \Delta) \stackrel{\mu^i}{\longrightarrow} H^i(\Gamma , n\Delta)
\end{equation}
and
\begin{equation}\label{HSE:2}
H^{i-1}(\Gamma , \Delta/n\Delta) \longrightarrow H^i(\Gamma , n\Delta) \stackrel{\nu^i}{\longrightarrow} H^i(\Gamma , \Delta).
\end{equation}
Since $\vert \Gamma \vert = n$ and $i \geq 1$, it is well known that the group $H^i(\Gamma , \Delta)$ is annihilated by $n$ (see, for example, \cite[Ch. II, Corollary 1.31]{MilneCFT}). In other words, we have
\begin{equation}\label{HSE:2a}
\nu^i \circ \mu^i = 0.
\end{equation}
It follows from (\ref{HSE:2}) and the finiteness of $\Delta/n\Delta$ that $\ker \nu^i$ is finite; hence, by (\ref{HSE:2a}),
$\im \mu^i$ is finite as well. On the other hand, we conclude from (\ref{HSE:1}) that $\im \mu^i \simeq {\mathrm {coker}}\: \lambda^i$. In view of the finiteness of ${}_n\Delta$, the group $H^i(\Gamma , {}_n\Delta)$ is finite, so the finiteness of ${\mathrm {coker}}\: \lambda^i$ readily yields the finiteness of $H^i(\Gamma , \Delta)$.
\end{proof}

\vskip2mm

\begin{remark}\label{R:HSL}
One point to note, which will be useful for us later, is that the assumptions of the lemma depend only on the structure of $\Delta$ as an {\it abelian group}, and not as a $\Gamma$-module.
\end{remark}

\vskip2mm

We now fix the notations that will be needed in the proof of Theorem \ref{T-HigherSha}.
Let $K = k(X)$ be the function field of a smooth geometrically integral variety $X$ over a field $k$ that has characteristic 0 and is of type (F) and set $V$ to be the set of geometric places of $K$ corresponding to the prime divisors of $X$. For our purposes, we may (and will) assume that $X$ is affine --- this can always be achieved by replacing $X$ with a suitable affine open $k$-subvariety, which will only shrink the corresponding set of geometric places of $K$.

We note that $V$ automatically satisfies condition (*) introduced in \S2.2, so, we are in a position to use the adelic set-up recalled there. In particular, if $L/K$ is a finite Galois extension, let $\bar{V}$ denote the set of all extensions to $L$ of places in $V$.
Then, for any $K$-torus $T$, the diagonal embedding $T(L) \hookrightarrow T(\mathbb{A}(L, \bar{V}))$
gives rise to group homomorphisms
$$
\nu^i_{L/K} \colon H^i(L/K, T) \to H^i(L/K, T(\mathbb{A}(L, \bar{V})))
$$
for all $i \geq 1.$

We first note that the argument used in the proof of \cite[Lemma 4.3]{RR-Tori} yields the following.

\begin{prop}\label{P:AdeleHigherSha}
For any $K$-torus $T$ and  any finite Galois extension $L/K$, we have
$$
\Sha^i(L/K, T, V) = \ker \left(H^i(L/K, T) \stackrel{\nu^i_{L/K}}{\longrightarrow} H^i(L/K, T(\mathbb{A}(L, \bar{V}))) \right)
$$
for all $i \geq 1.$
\end{prop}

\vskip2mm

\noindent {\it Proof of Theorem \ref{T-HigherSha}.} Let $T$ be a torus over $K = k(X)$ and let $L/K$ be a finite Galois extension that splits $T$, so that
$$
T \times_K L \simeq \mathbb{G}_{m,L}^d, \ \ \ d = \dim T.
$$
As above, let $\bar{V}$ denote the set of all extensions to $L$ of places in $V$. 
It is easy to see that in geometric terms, the elements of $\bar{V}$ correspond to the prime divisors of the normalization $\tilde{X}$ of $X$ in $L$. We note that $\tilde{X}$ is defined over the algebraic closure $\ell$ of $k$ in $L$ and we have $L = \ell(\tilde{X}).$ Since $\ell/k$ is a finite extension and $k$ is of type (F), it follows that $\ell$ is also of type (F).
\begin{lemma}\label{L:PicFinite}
Set $\Delta = \mathrm{Pic}(L,\bar{V}).$ Then for any integer $n \geq 1$, the groups ${}_n\Delta$ and $\Delta/n \Delta$ are finite.
\end{lemma}
\begin{proof}
Let $\tilde{X}_{\rm sing}$ denote the singular locus of $\tilde{X}$ and set $U = \tilde{X} \setminus \tilde{X}_{\rm sing}.$ Since $\tilde{X}$ is normal, we have $\mathrm{codim}_{\tilde{X}} \tilde{X}_{\rm sing} \geq 2$ by Serre's normality criterion (cf. \cite[Proposition B.78]{GW}), and hence $\mathrm{Pic}(L,\bar{V}) = \mathrm{Pic}(L,\bar{V}_U)$, where $\bar{V}_U$ is the set of discrete valuations of $L$ corresponding to the prime divisors of $U$. On the other hand, due to the fact that $U$ is smooth, we can identify $\mathrm{Pic}~(L,\bar{V}_U)$
with the usual Picard group $\mathrm{Pic}~U$ of $U$. Now, since
${\rm char}~\ell = 0$ and $\ell$ is of type (F) (as noted above), according to Theorem \ref{T-EtaleFinite},
the \'etale cohomology groups $\he^i(U, \mu_{n, U})$ are finite for all $i \geq 0$ and all $n \geq 1$. One then utilizes the Kummer sequence
$$
1 \to \mu_{n, U} \to \mathbb{G}_{m,U} \stackrel{[n]}{\longrightarrow} \mathbb{G}_{m, U} \to 1.
$$
to obtain the finiteness of ${}_n\Delta$ and $\Delta/n \Delta$ as in the proof of Theorem \ref{T-Surface} (see \S\ref{S:Surface}).
\end{proof}

Let us now consider the following short exact sequence of $\Ga(L/K)$-modules
$$
1 \to T(\mathbb{A}^{\infty}(L, \bar{V})) T(L) \longrightarrow T(\mathbb{A}(L, \bar{V})) \longrightarrow T(\mathbb{A}(L, \bar{V})) / (T(\mathbb{A}^{\infty}(L, \bar{V})) T(L)) \to 1
$$
and the corresponding exact sequence in cohomology
\begin{multline*}
H^{i-1}\big(L/K,\ T(\mathbb{A}(L,\bar V))\,/\,T(\mathbb{A}^{\infty}(L,\bar V))\,T(L)\big) \\
\longrightarrow\ H^i\big(L/K,\ T(\mathbb{A}^{\infty}(L,\bar V))\,T(L)\big)\ \stackrel{\alpha^i}{\longrightarrow}\ H^i\big(L/K,\ T(\mathbb{A}(L,\bar V))\big).
\end{multline*}
In view of the identification
$$
T(\mathbb{A}(L,\bar{V})) / T(\mathbb{A}^{\infty}(L, \bar{V})) T(L) \simeq \left( \mathbb{I}(L, \bar{V}) / \mathbb{I}^{\infty}(L, \bar{V}) L^{\times} \right)^d \simeq (\mathrm{Pic}(L,\bar{V}))^d,
$$
recalled in (\ref{E:PicIdeles}) in \S2.3, it follows from Lemmas \ref{HSL:1} and \ref{L:PicFinite} (see also Remark \ref{R:HSL}) that the group $$H^{i-1}(L/K , T(\mathbb{A}(L, \bar{V})) / (T(\mathbb{A}^{\infty}(L, \bar{V})) T(L)))$$ is finite, and hence
$\ker \alpha^i$ is finite for all $i \geq 2.$ On the other hand, since
$$
\Sha^i(L/K, T, V) = \ker \left(H^i(L/K, T) \stackrel{\nu^i_{L/K}}{\longrightarrow} H^i(L/K, T(\mathbb{A}(L, \bar{V}))) \right)
$$
by Proposition \ref{P:AdeleHigherSha}, we see that to complete the proof of Theorem \ref{T-HigherSha}, it suffices to establish the finiteness of
$$
\ker \left(\tau^i_{L/K} \colon H^i(L/K , T) \longrightarrow H^i(L/K , T(\mathbb{A}^{\infty}(L, \bar{V})) T(L)) \right),
$$
where $\tau^i_{L/K}$ is induced by the natural embedding
$$
T(L) \hookrightarrow T(\mathbb{A}^{\infty}(L, \bar{V})) T(L).
$$
For this, we let $R = T(L) \cap T(\mathbb{A}^{\infty}(L, \bar{V}))$ and consider the exact sequence
$$
1 \to R \longrightarrow T(L) \times T(\mathbb{A}^{\infty}(L, \bar{V})) \stackrel{\varphi \cdot \psi}{\longrightarrow} T(L) T(\mathbb{A}^{\infty}(L, \bar{V})) \to 1,
$$
where the first map is given by $r \mapsto (r , r^{-1})$ and the second is the product of the natural embeddings
$$
\varphi \colon T(L) \hookrightarrow T(L) T(\mathbb{A}^{\infty}(L, \bar{V})) \ \ \text{and} \ \ \psi \colon T(\mathbb{A}^{\infty}(L, \bar{V})) \hookrightarrow T(L) T(\mathbb{A}^{\infty}(L, \bar{V})).
$$
We then have the corresponding cohomological exact sequence
\begin{equation}\label{E-HS5}
H^i(L/K , R) \longrightarrow H^i(L/K , T) \times H^i\big(L/K , T(\mathbb{A}^{\infty}(L, \bar{V}))\big) \stackrel{\varphi^i + \psi^i}{\longrightarrow}
H^i\big(L/K , T(L) T(\mathbb{A}^{\infty}(L, \bar{V}))\big),
\end{equation}
where
$$
\varphi^i \colon H^i(L/K , T) \longrightarrow H^i(L/K , T(L) T(\mathbb{A}^{\infty}(L, \bar{V})))
$$
and
$$
\psi^i \colon H^i(L/K , T(\mathbb{A}^{\infty}(L, \bar{V}))) \longrightarrow H^i(L/K , T(L) T(\mathbb{A}^{\infty}(L, \bar{V})))
$$
are induced by $\varphi$ and $\psi$, respectively. Clearly, $\ker \tau^i_{L/K}$ embeds into $\ker (\varphi^i + \psi^i)$ via $x \mapsto (x , 0)$, so it is enough to establish the finiteness of $\ker (\varphi^i + \psi^i)$. In view of (\ref{E-HS5}), it suffices to show that $H^i(L/K, R)$ is finite, which we will do using
Lemma \ref{HSL:1}. To verify the assumptions of the lemma, we first note that, as an abelian group, $R$ is isomorphic to $(L^{\times} \cap \mathbb{I}^{\infty}(L, \bar{V}))^d$. Thus, we need the following.

\begin{lemma}\label{HSL:3}
Set $\Theta = L^{\times} \cap \mathbb{I}^{\infty}(L, \bar{V})$. Then for any integer $n \geq 1$, the groups ${}_n\Theta$ and $\Theta/ n\Theta$ are finite.

\end{lemma}
\begin{proof}
As noted above, we have
$L = \ell(\tilde{X})$, where $\tilde{X}$ is the normalization of $X$ in $L$ and $\ell$ is the algebraic closure of $k$ in $L$.
By our assumption, $X$ is affine, so $\tilde{X}$ is also affine, and hence we can embed it as a dense open subset into a projective $\ell$-variety $Y$. Furthermore, since $\tilde{X}$ is normal, and it is well-known that the normalization of a projective variety is projective (see, for example, \cite[Tag 0BXQ, Lemma 33.27.2]{Stacks}), we can assume that $Y$ is normal as well. We have $\ell(Y) = \ell(\tilde{X}) = L$, and we denote by $\bar{V}_Y$ the set of discrete valuations of $L$ associated with the prime divisors on $Y$. Then $\bar{V}_Y \setminus \bar{V}$ consists of finitely many valuations $v_1, \dots, v_r$ corresponding to the prime divisors on $Y$ not contained in $\tilde{X}.$ Consider the map
$$
\varkappa \colon \Theta \longrightarrow \mathbb{Z}^r, \ \ \ x \mapsto (v_1(x), \ldots , v_r(x)).
$$
Since $Y$ is normal and complete, any function $a \in L$ that satisfies $v(a) \geq 0$ for all $v \in \bar{V}_Y$ belongs to $\ell$ (see, for example, \cite[Theorem 6.45]{GW} or \cite[Proposition 6.3A]{Hartshorne}). It follows that
$$\ker \varkappa = \{ a \in L^{\times} \ \vert \ v(a) = 0 \ \ \text{for all} \ \ v \in \bar{V}_Y\}$$
reduces to $\ell^{\times}$. On the other hand, $\im \varkappa$ is isomorphic to $\mathbb{Z}^s$ for some $s \leq r$. Thus,
$$
\Theta \simeq \ell^{\times} \times \mathbb{Z}^s.
$$
Clearly, ${}_n\ell^{\times}$ is finite. Furthermore, since $\ell$ is a finite extension of $k$, it is also of type (F),
and
therefore $\ell^{\times}/ (\ell^{\times})^n$ is finite (see \cite[Lemma 2.4]{IRap}). We conclude that ${}_n\Theta$ and $\Theta/ n\Theta$ are both finite, as needed.
\end{proof}

Thus, $\Theta$ satisfies the assumptions of Lemma \ref{HSL:1}, and hence the groups $H^i(L/K, R)$ are finite for all $i \geq 1.$ As noted above, this completes the proof of Theorem \ref{T-HigherSha}. \hfill $\Box$

\vskip4mm

\noindent {\bf Acknowledgements.} I was partially supported by NSF grant DMS-2154408 during the preparation of this paper.

\bibliographystyle{amsplain}

\end{document}